\documentclass[11pt,a4paper]{article}
\usepackage{a4wide}
\usepackage[T1]{fontenc}
\usepackage{latexsym}
\usepackage{amsmath}
\usepackage{amsthm}
\usepackage{amssymb,enumerate}
\usepackage[usenames]{color}
\usepackage{graphicx}
\usepackage{fancybox}
\usepackage[normalem]{ulem}
\usepackage{hyperref}
\theoremstyle{plain}
\newtheorem{theo}{Theorem}[section]
\newtheorem{lemma}[theo]{Lemma}
\newtheorem{prop}[theo]{Proposition}
\newtheorem{coro}[theo]{Corollary}
\newtheorem*{coro*}{Corollary}
\theoremstyle{definition}
\newtheorem{defi}[theo]{Definition}
\newtheorem{rema}[theo]{Remark}

\newenvironment{proof1}{\medskip\par\noindent{\bf Proof.}}{\hfill $\Box$
\medskip\par}

\newcommand{\RR}{\mathbb{R}}

\newcommand{\N}{\mathbb{N}}
\newcommand{\R}{\mathbb{R}}
\newcommand{\C}{\mathbb{C}}

\newcommand{\ka}{\kappa}
\newcommand{\al}{\alpha}
\newcommand{\be}{\beta}

\newcommand{\bM}{\mathbb{M}}
\newcommand{\M}{\mathbb{M}}
\newcommand{\hM}{\widehat{\M}}
\newcommand{\bL}{\mathbb{L}}
\newcommand{\m}{\boldsymbol{m}}
\newcommand{\bl}{\boldsymbol{\ell}}
\newcommand{\mo}{\mathfrak{m}}
\newcommand{\flo}[1]{\left\lfloor {#1} \right\rfloor}

\def\a{\alpha}

\def\o{\omega}

\def\ga{\gamma}
\def\L{\mathbb{L}}
\def\bm{\boldsymbol{m}}
\def\Ms{{\mathbb{M}}_{q,\sigma}}

\definecolor{azulosc}{rgb}{0.2,0.1,0.7}
\definecolor{granate}{rgb}{0.6,0,0.3} 

\definecolor{verdeosc}{RGB}{46, 139, 87} 
\definecolor{rojo}{RGB}{219,0,0}                

\makeatletter
\DeclareRobustCommand\widecheck[1]{{\mathpalette\@widecheck{#1}}}
\def\@widecheck#1#2{%
    \setbox\z@\hbox{\m@th$#1#2$}%
    \setbox\tw@\hbox{\m@th$#1%
       \widehat{%
          \vrule\@width\z@\@height\ht\z@
          \vrule\@height\z@\@width\wd\z@}$}%
    \dp\tw@-\ht\z@
    \@tempdima\ht\z@ \advance\@tempdima2\ht\tw@ \divide\@tempdima\thr@@
    \setbox\tw@\hbox{%
       \raise\@tempdima\hbox{\scalebox{1}[-1]{\lower\@tempdima\box
\tw@}}}%
    {\ooalign{\box\tw@ \cr \box\z@}}}
\makeatother

\begin{document}

\title{Optimal flat functions in Carleman-Roumieu ultraholomorphic classes in sectors}
\author{Javier Jim\'enez-Garrido \and Ignacio Miguel-Cantero \and Javier Sanz \and Gerhard Schindl}
\date{\today}

\maketitle

\begin{abstract}
We construct optimal flat functions in Carleman-Roumieu ultraholomorphic classes associated to general strongly nonquasianalytic weight sequences, and defined on sectors of suitably restricted opening. A general procedure is presented in order to obtain linear continuous extension operators, right inverses of the Borel map, for the case of regular weight sequences in the sense of Dyn'kin. Finally, we discuss some examples (including the well-known $q$-Gevrey case) where such optimal flat functions can be obtained in a more explicit way.
\par\medskip

\noindent Key words: Carleman-Roumieu ultraholomorphic classes, asymptotic expansions, linear extension operators.
\par
\medskip
\noindent 2020 MSC: Primary 47A57; secondary 44A10, 46E10. \end{abstract}

\section{Introduction}
The asymptotic Borel map sends a function, admitting an asymptotic expansion in a sectorial region, into the formal power series providing such expansion. In many instances it is important to decide about the injectivity and surjectivity of this map when considered between so-called Carleman-Roumieu ultraholomorphic classes and the corresponding class of formal series, defined by restricting the growth of some of the characteristic data of their elements (the derivatives of the functions, the remainders in the expansion, or the coefficients of the series) in terms of a given weight sequence $\M=(M_p)_{p\in\N_0}$ of positive real numbers (see Subsection~\ref{subsectCarlemanclasses} for the definition of such classes). While the injectivity has been fully characterized for sectorial regions and general weight sequences~\cite{Mandelbrojt,Salinas,JimenezSanzSchindlInjectSurject}, the surjectivity problem is still under study.
The classical Borel-Ritt-Gevrey theorem of B.~Malgrange and J.-P.~Ramis~\cite{Ramis1}, solving the case of Gevrey asymptotics (for which $\M=(p!^{\alpha})_{p\in\N_0}$, $\alpha>0$), was partially extended to different more general situations by J.~Schmets and M.~Valdivia~\cite{SchmetsValdivia00}, V.~Thilliez~\cite{Thilliez95,Thilliez03} and the authors~\cite{SanzFlatProxOrder,JimenezSanzSchindlInjectSurject,JimenezSanzSchindlSurjectDC}. Summing up, it is known now that the strong nonquasianalyticity condition (snq) for $\M$, equivalent to the fact that the index $\gamma(\M)$ introduced by V. Thilliez is positive (see Subsection~\ref{subsectIndexGammaM}), is indeed necessary for surjectivity. Moreover, for an unbounded sector $S_{\gamma}$ of opening $\pi\gamma$ ($\gamma>0$) in the Riemann surface of the logarithm and for regular weight sequences in the sense of E.~M.~Dyn'kin~\cite{Dynkin80} (see Subsection~\ref{subsectstrregseq} for the definitions),
the Borel map is surjective whenever $\gamma<\gamma(\M)$, while it is not for $\gamma>\gamma(\M)$ (the situation for $\gamma=\gamma(\M)$ is still unclear in general). It is important to note that
the current proof of surjectivity in this situation is not constructive, but rests on the characterization, by abstract functional-analytic techniques, of the surjectivity of the Stieltjes moment mapping in Gelfand-Shilov spaces defined by regular sequences due to A. Debrouwere~\cite{momentsdebrouwere}. This information is transferred into the asymptotic framework in a halfplane by means of the Fourier transform, and in~\cite{JimenezSanzSchindlSurjectDC} Laplace and Borel transforms of arbitrary order allow to conclude for general sectors.
However, in the particular case of classes given by strongly regular sequences in the sense of V. Thilliez, the proof of surjectivity of the Borel map~\cite{Thilliez03} rests on the construction of optimal flat functions in suitable sectors and a double application of Whitney extension results. Subsequently, A. Lastra, S. Malek and the third author~\cite{LastraMalekSanzContinuousRightLaplace} reproved surjectivity in a more explicit way by means of formal Borel- and truncated Laplace-like transforms, defined from suitable kernel functions obtained from those optimal flat functions.

The first aim of this paper is to construct such optimal flat functions for Carleman-Roumieu ultraholomorphic classes defined by general weight sequences (not just strongly regular ones) and in sectors $S_\gamma$ with $\gamma<\gamma(\M)$. The key idea comes from a recent preprint by D.~N.~ Nenning, A.~Rainer and the fourth author~\cite{NenningRainerSchindl_preprint}, where they have studied the mixed Borel problem in Beurling ultradifferentiable classes. They consider a mixed condition inspired by a related one (see~\eqref{eq.Langenbruch} in this paper) appearing in a work of M.~Langenbruch~\cite{Langenbruch94}. It turns out that the condition of Langenbruch is, under natural hypotheses, equivalent to the fact that  $\gamma(\M)>1$, and it is crucial in order to construct optimal flat functions in a halfplane by means of the classical harmonic extension of the associated function $\omega_{\M}$. A ramification process provides then optimal flat functions in the general situation.
These results completely close the problem of the explicit construction of optimal flat functions in sectors of appropriate opening for classes defined in terms of a general weight sequence. Moreover, the constructive techniques developed in this paper could be used in other contexts where weighted structures appear.

Secondly, for ultraholomorphic classes defined by regular sequences we establish the connection with the surjectivity of the Borel map by providing a constructive technique for the corresponding extension results, in the same vein as in~\cite{LastraMalekSanzContinuousRightLaplace}. For sake of completeness, in the case of strongly regular sequences we also give an alternative approach, based on the work of J.~Bruna~\cite{Brunaext80}, to this connection.

In order to highlight the power of the technique in concrete situations, we will also present a family of (non strongly) regular sequences for which such optimal flat functions can be provided in any sector of the Riemann surface of the logarithm (what agrees with the fact that the index $\gamma(\M)$ is in this case equal to $\infty$), resting on precise estimates for the associated function $\omega_{\M}$ instead of appealing to its harmonic extension.
We end by showing how optimal flat functions and extension results can be obtained for convolved sequences, in case the factor sequences admit such constructions separately. Some examples are commented on in regard with this technique.

The paper is organized as follows. Section~\ref{sectPrelimin} consists of all the preliminary information concerning weight sequences and some indices or auxiliary functions associated with them, and the main facts about ultraholomorphic classes and the (asymptotic) Borel map defined for them. In Section~\ref{sectFlatFunctions} we define optimal flat functions and carefully detail their construction for general weight sequences. Next, we show that their existence entails the
surjectivity of the Borel map in ultraholomorphic classes defined by regular sequences. In the particular case of sequences of moderate growth, different statements are presented relating the property of strong non-quasianalyticity to the existence of such flat functions. In Section~\ref{sectConstrOptFlatNon_SR_Seq} we give a family of sequences (among which the classical $q$-Gevrey sequences are found) for which optimal flat functions can be constructed in a more explicit way.
We need to work first in $\C\setminus(-\infty,0]$, and then apply a ramification in order to reason for arbitrary sectors in the Riemann surface of the logarithm.
Finally, the last section is devoted to the work with convolved sequences.

\section{Preliminaries}\label{sectPrelimin}

\subsection{Weight sequences and their properties}\label{subsectstrregseq}

We set $\N:=\{1,2,...\}$, $\N_{0}:=\N\cup\{0\}$.
In what follows, $\bM=(M_p)_{p\in\N_0}$ will always stand for a sequence of positive real numbers with $M_0=1$. We define its {\it sequence of quotients}  $\m=(m_p)_{p\in\N_0}$ by
$m_p:=M_{p+1}/M_p$, $p\in \N_0$; the knowledge of $\M$ amounts to that of $\m$, since $M_p=m_0\cdots m_{p-1}$, $p\in\N$.
The following properties for a sequence will play a role in this paper:\par
(i) $\M$ is \emph{logarithmically convex} (for short, (lc)) if
$M_{p}^{2}\le M_{p-1}M_{p+1}$, $p\in\N$.\par
(ii) $\M$ is \emph{stable under differential operators} or satisfies the \emph{derivation closedness condition} (briefly, (dc)) if there exists $D>0$ such that $M_{p+1}\leq D^{p+1} M_{p}$, $p\in\N_{0}$.\par
(iii) $\M$ is of, or has, \emph{moderate growth} (for the sake of brevity, (mg)) if there exists $A>0$ such that
$M_{p+q}\le A^{p+q}M_{p}M_{q}$, $p,q\in\N_0$.\par
(iv) $\M$ satisfies the condition (nq) of \emph{non-quasianalyticity} if
$$\sum_{p=0}^{\infty}\frac{M_{p}}{(p+1)M_{p+1}}<+\infty.
$$
\indent (v) Finally, $\M$ satisfies the condition (snq) of \emph{strong non-quasianalyticity} if there exists $B>0$ such that
$$
\sum^\infty_{q= p}\frac{M_{q}}{(q+1)M_{q+1}}\le B\frac{M_{p}}{M_{p+1}},\qquad p\in\N_0.$$

It is convenient to introduce the notation $\hM:=(p!M_p)_{p\in\N_0}$.
All these properties are preserved when passing from $\M$ to $\hM$.
In the classical work of H.~Komatsu~\cite{komatsu}, the properties (lc), (dc) and (mg) are denoted by $(M.1)$, $(M.2)'$ and $(M.2)$, respectively, while
(nq) and (snq) for $\M$ are the same as properties $(M.3)'$ and $(M.3)$ for $\widehat{\M}$, respectively.
Obviously, (mg) implies (dc).

The sequence of quotients $\m$ is nondecreasing if and only if $\M$ is (lc). In this case, it is well-known that $(M_p)^{1/p}\leq m_{p-1}$ for every $p\in\N$, the sequence $((M_p)^{1/p})_{p\in\N}$ is nondecreasing, and $\lim_{p\to\infty} (M_p)^{1/p}= \infty$ if and only if $\lim_{p\to\infty} m_p= \infty$. In order to avoid trivial situations, we will restrict from now on to
(lc) sequences $\M$ such that $\lim_{p\to\infty} m_p =\infty$, which will be called \emph{weight sequences}.

Following E.~M.~Dyn'kin~\cite{Dynkin80}, if $\M$ is a weight sequence and satisfies (dc), we say $\hM$ is \emph{regular}. According to V.~Thilliez~\cite{Thilliez03}, if $\M$ satisfies (lc), (mg) and (snq), we say $\M$ is \emph{strongly regular}; in this case $\M$ is a weight sequence, and the corresponding $\hM$ is regular.

We mention some interesting examples. In particular, those in (i) and (iii) appear in the applications of summability theory to the study of formal power series solutions for different kinds of equations.
\begin{enumerate}[(i)]
\item The sequences $\M_{\al,\be}:=\big(p!^{\al}\prod_{m=0}^p\log^{\be}(e+m)\big)_{p\in\N_0}$, where $\al>0$ and $\be\in\R$, are strongly regular (in case $\be<0$, the first terms of the sequence have to be suitably modified in order to ensure (lc)). In case $\be=0$, we have the best known example of a strongly regular sequence, $\M_{\al}:=\M_{\al,0}=(p!^{\al})_{p\in\N_{0}}$, called the \emph{Gevrey sequence of order $\al$}.

\item The sequence $\M_{0,\be}:=(\prod_{m=0}^p\log^{\be}(e+m))_{p\in\N_0}$, with $\be>0$, satisfies (lc) and (mg), and $\m$ tends to infinity, but (snq) is not satisfied.
\item For $q>1$ and $1<\sigma\le 2$, $\M_{q,\sigma}:=(q^{p^\sigma})_{p\in\N_0}$ satisfies (lc), (dc) and (snq), but not (mg). In case $\sigma=2$, we get the well-known \emph{$q$-Gevrey sequence}.
\end{enumerate}

Two sequences $\bM=(M_{p})_{p\in\N_0}$ and $\bL=(L_{p})_{p\in\N_0}$ of positive real numbers,
with $M_0=L_0=1$ and with respective quotient sequences $\m$ and $\bl$, are said to be \emph{equivalent}, and we write $\M\approx\bL$, if there exist positive constants $A,B$ such that $A^pM_p\le L_p\le B^pM_p$, $p\in\N_0$. They are said to be \emph{strongly equivalent}, denoted by $\m\simeq\bl$, if there exist positive constants $a,b$ such that
$am_p\le \ell_p\le bm_p$, $p\in\N_0$. Whenever $\m\simeq\bl$ we have $\M\approx\bL$ but, in general, not conversely.

In case $M_0$ or $L_0$ is not equal to 1, the previous definitions of equivalence are meant to deal with the normalized sequences $(M_{p}/M_0)_{p\in\N_0}$ or $(L_{p}/L_0)_{p\in\N_0}$.

\subsection{Index $\ga(\M)$ and auxiliary functions for weight sequences and functions}\label{subsectIndexGammaM}

The index $\gamma(\M)$, introduced by V.~Thilliez~\cite[Sect.\ 1.3]{Thilliez03} for strongly regular sequences $\M$, can be equally defined for (lc) sequences, and it may be equivalently expressed by different conditions:
\begin{enumerate}[(i)]
\item A sequence $(c_p)_{p\in\N_0}$ is \emph{almost increasing} if there exists $a>0$ such that for every $p\in\N_0$ we have that $c_p\leq a c_q $ for every $ q\geq p$.
It was proved in~\cite{JimenezSanzSRSPO,JimenezSanzSchindlIndex} that for any weight sequence $\M$ one has
\begin{equation}\label{equa.indice.gammaM.casicrec}
\gamma(\bM)=\sup\{\gamma>0:(m_{p}/(p+1)^\gamma)_{p\in\N_0}\hbox{ is almost increasing} \}\in[0,\infty].
\end{equation}
\item For any $\be>0$ we say that $\m$ satisfies the condition $(\gamma_{\be})$ if there exists $A>0$ such that
\begin{equation*}
\sum^\infty_{\ell=p} \frac{1}{(m_\ell)^{1/\be}}\leq \frac{A (p+1) }{(m_p)^{1/\be}},  \qquad p\in\N_0.\tag{$\gamma_\be$}
\end{equation*}
In~\cite{PhDJimenez,JimenezSanzSchindlIndex} it is proved that for a weight sequence $\M$,
\begin{equation}\label{equa.indice.gammaM.gamma_r}
\gamma(\M)=\sup\{\be>0; \,\, \m \,\, \text{satisfies } (\gamma_{\be})\,\};\quad \gamma(\M)>\be\iff\m\text{ satisfies }(\gamma_\be).
\end{equation}
\end{enumerate}

If we observe that the condition (snq) for $\M$ is precisely $(\gamma_1)$ for $\widehat{m}$, the sequence of quotients for $\hM$, and that $\ga(\hM)=\ga(\M)+1$ (this is clear from~\eqref{equa.indice.gammaM.casicrec}), we deduce from the second statement in~\eqref{equa.indice.gammaM.gamma_r} that
\begin{equation}\label{eq.snqiffgammaMpositive}
\textrm{$\M$ satisfies (snq) if, and only if, }\ga(\M)>0.
\end{equation}

Given a weight sequence $\M=(M_p)_{p\in\N_0}$, we write
$$\o_{\M}(t):=\sup_{p\in\N_0} \ln\left(\frac{t^{p}}{M_{p}}\right),\qquad t>0,$$
and $\o_{\M}(0)=0$. This is the classical
(continuous, nondecreasing)
function  associated with the sequence $\M$, see~\cite{komatsu}.

Another associated function will play a key-role, namely
\begin{equation*}
h_\M(t):=\inf_{p\in\N_0}M_p t^p, \quad t>0.
\end{equation*}
The functions $h_\M$ and $\omega_\M$ are related by
\begin{equation}\label{functionhequ2}
h_\M(t)=\exp(-\omega_\M(1/t)),\quad t>0.
\end{equation}

In~\cite[Prop.~3.2]{komatsu} we find that, for a weight sequence $\M$,
\begin{equation}\label{eq.MpfromomegaM}
M_p=\sup_{t>0}t^p\exp(-\o_{\M}(t))=\sup_{t>0}t^p h_{\M}(1/t),\quad p\in\N_0.
\end{equation}

We record for the future some elementary facts about $h_{\M}$.

\begin{lemma}\label{functionhproperties}
Let $\M=(M_p)_{p\in\N_0}$ be a weight sequence, then:
\begin{itemize}
\item[$(i)$] $t\in(0,\infty)\mapsto h_\M(t)$ is nondecreasing and continuous,
\item[$(ii)$] $h_{\M}(t)\le 1$ for all $t>0$, $h_\M(t)=1$ for all $t$ sufficiently large and $\lim_{t\rightarrow 0}h_\M(t)=0$.
\end{itemize}
\end{lemma}

We also introduce the counting function $\nu_{\m}$ for the sequence $\m$,
\begin{equation}\label{countinfunctionnu}
\nu_{\m}(\lambda):=\#\{p\in\N_0: m_p\le\lambda\}.
\end{equation}
If $\M$ is a weight sequence, then the functions $\nu_{\m}$ and $\omega_{\M}$ are related by the following useful integral representation formula, e.g. see~\cite{Mandelbrojt}
and~\cite[$(3.11)$]{komatsu}:
\begin{equation}\label{assofuncintegral}
\omega_\M(x)= \int_0^{x}\frac{\nu_{\m}(\lambda)}{\lambda}d\lambda= \int_{m_0}^{x}\frac{\nu_{\m}(\lambda)}{\lambda}d\lambda,\quad x>0.
\end{equation}

In \cite{JimenezSanzSchindlIndex}, the nature of the index $\gamma(\M)$, fundamental in the study of the surjectivity of the asymptotic Borel map, is explained. More precisely, it is shown that $\gamma(\M)$ is the lower Matuszewska index  of $\m$. In addition, the relation between ORV-indices of $\m$ and $\nu_{\m}$ is clarified and from this connection we characterized some properties of $\nu_{\m}$ that will be important for our aim.

\begin{lemma}\label{functionNUproperties}
Let $\M=(M_p)_{p\in\N_0}$ be a weight sequence, then:
\begin{itemize}
\item[$(i)$] $\gamma(\M)>0$ if and only if $\nu_{\m}$ satisfies the condition $\nu_{\m}(2t)=O(\nu_{\m}(t))$ as $t$ tends to $\infty$.
\item[$(ii)$] $\gamma(\M)>1$ if and only if $\nu_{\m}$ satisfies the condition $(\omega_{\text{snq}})$, i.e., there exists $D>0$ such that
	\begin{equation*}
		\int_1^\infty \frac{\nu_{\m}(ys)}{s^2}\,ds\leq D\nu_{\m}(y)+D, \qquad y\geq 0.
	\end{equation*}
\end{itemize}
\end{lemma}

\begin{proof} $(i)$ follows by \eqref{eq.snqiffgammaMpositive} and~\cite[Coro.\ 4.2.(ii)]{JimenezSanzSchindlIndex}. $(ii)$  holds true by combining ~\cite[Lemma \ 2.10]{JimenezSanzSchindlIndex}, ~\cite[Coro.\ 2.13]{JimenezSanzSchindlIndex}, ~\cite[Thm. 3.10]{JimenezSanzSchindlIndex} and~\cite[Prop.\ 4.1]{JimenezSanzSchindlIndex}.
\end{proof}

We conclude this subsection by introducing the harmonic extension and a particular majorant of a nondecreasing nonquasianalytic function.

A nondecreasing (or even just measurable) function $\sigma:[0,\infty)\to [0,\infty)$ satisfies the \emph{nonquasianalyticity} property $(\omega_{nq})$ (and we say $\sigma$ is \emph{nonquasianalytic}) if
	$$\int_{1}^{\infty}\frac{\sigma(t)}{t^2}\,dt<\infty.$$

Let $\sigma:[0,\infty)\to [0,\infty)$ be a nondecreasing nonquasianalytic function. The \emph{harmonic extension} $P_\sigma$ of $\sigma$ to the open upper and lower halfplanes of $\C$ is defined by
	\begin{equation}\label{harmonicextension}
		P_\sigma(x+iy)=\begin{cases}
			\quad\sigma(|x|)&\text{if } x\in\R,\,y=0,\\
\displaystyle\frac{|y|}{\pi}\int_{-\infty}^{\infty}\frac{\sigma(|t|)}{(t-x)^2+y^2}dt& \text{if }x\in\R,\,y\neq0.
		\end{cases}
	\end{equation}

For every $z\in\C$ one has (see, for example,~\cite[Remark 3.2]{BonetMeiseTaylor92} or~\cite[Prop.\ 5.5]{NenningRainerSchindl_preprint}):
	\begin{equation}\label{omegapbound}
		\sigma(|z|)\leq P_{\sigma}(z).
	\end{equation}

We list some basic properties of the harmonic extension that will be used later:
\begin{enumerate}
	\item[(1)] $\sigma_1 \leq \sigma_2$ implies $P_{\sigma_1}\leq P_{\sigma_2}$.
	\item[(2)] $\lambda P_{\sigma_1}(z)+\mu P_{\sigma_2}(z)=P_{\lambda\sigma_1+\mu\sigma_2}(z)$, $\lambda,\mu \in\RR$.
	\item[(3)] $P_{t\mapsto\sigma(Ct)}(z)=P_{\sigma}(Cz)$, $C>0$.
\end{enumerate}

Another important auxiliary function appears in the study of extension results in Braun-Meise-Taylor ultradifferentiable classes, defined in terms of weight functions (see, for example,~\cite{MeiseTaylor88,BonetMeiseTaylor92} and the references therein).

Let $\sigma:[0,\infty)\to [0,\infty)$ be a nondecreasing and nonquasianalytic function. Then, the function $\ka_\sigma$ is defined by
	$$\ka_\sigma(y)=\int_1^\infty \frac{\sigma(ys)}{s^2}\,ds,\qquad y\geq0,$$
and satisfies
\begin{equation}\label{sigmakappabound}
\sigma(y)\leq \ka_\sigma(y),\quad y\ge 0.
\end{equation}
If $\sigma$ is also continuous, then $\ka_\sigma$ is concave, cf. the proof of $(3)\Rightarrow(4)$ in Prop. 1.3 in \cite{MeiseTaylor88}.

In particular, consider a weight sequence $\M$ such that $\sum_{p=0}^\infty 1/m_p<\infty$ (this is condition $(M3)'$ in~\cite{komatsu}); in other words, the sequence $\widecheck{\M}:=(M_p/p!)_{p\in\N_0}$ satisfies (nq). According to~\cite[Lemma 4.1]{komatsu}, this property amounts to $\nu_m$ and/or $\omega_{\M}$ being nonquasianalytic. So, it makes sense to consider the concave function $\ka_{\omega_\M}$ associated with $\omega_{\M}$, and  $\ka_{\nu_{\m}}$ associated with the counting function $\nu_{\m}$. The equality
\begin{equation}\label{kappaomeganurelation}
		\ka_{\omega_\M}(y)=\omega_{\M}(y)+\ka_{\nu_{\m}}(y), \qquad y\geq 0,
\end{equation}
can be found on p. 58 in the proof of ~\cite[Prop.\ 4.4]{komatsu}.


%

\subsection{Asymptotic expansions, ultraholomorphic classes and the asymptotic Borel map}\label{subsectCarlemanclasses}

$\mathcal{R}$ stands for the Riemann surface of the logarithm.
$\C[[z]]$ is the space of formal power series in $z$ with complex coefficients.

For $\gamma>0$, we consider unbounded sectors bisected by direction 0,
$$S_{\gamma}:=\{z\in\mathcal{R}:|\hbox{arg}(z)|<\frac{\gamma\,\pi}{2}\}$$
or, in general, unbounded sectors with bisecting direction $d\in\R$ and opening $\ga\,\pi$,
$$S(d,\ga):=\{z\in\mathcal{R}:|\hbox{arg}(z)-d|<\frac{\ga\,\pi}{2}\}.$$

A sector $T$ is said to be a \emph{proper subsector} of a sector $S$ if $\overline{T}\subset S$ (where the closure of $T$ is taken in $\mathcal{R}$, and so the vertex of the sector is not under consideration).

In this paragraph $S$ is an unbounded sector and $\M$ a sequence. We start by recalling the concept of uniform asymptotic expansion.

We say a holomorphic function $f\colon S\to\C$ admits $\widehat{f}=\sum_{n\ge 0}a_nz^n\in\C[[z]]$ as its \emph{uniform $\M$-asymptotic expansion in $S$ (of type $1/A$ for some $A>0$)} if there exists $C>0$ such that for every $p\in\N_0$, one has
\begin{equation}\left|f(z)-\sum_{n=0}^{p-1}a_nz^n \right|\le CA^pM_{p}|z|^p,\qquad z\in S.\label{desarasintunifo}
\end{equation}
In this case we write $f\sim_{\M,A}^u\widehat{f}$ in $S$, and $\widetilde{\mathcal{A}}^u_{\M,A}(S)$ denotes the space of functions admitting uniform $\M$-asymptotic expansion of type $1/A$ in $S$, endowed with the norm
\begin{equation}\label{eq.NormUniformAsymptFixedType}
\left\|f\right\|_{\M,A,\overset{\sim}{u}}:=\sup_{z\in S,p\in\N_{0}}\frac{|f(z)-\sum_{k=0}^{p-1}a_kz^k|}{A^{p}M_{p}|z|^p},
\end{equation}
which makes it a Banach space. $\widetilde{\mathcal{A}}^u_{\{\M\}}(S)$ stands for the $(LB)$ space of functions admitting a uniform $\{\M\}$-asymptotic expansion in $S$, obtained as the union of the previous classes when $A$ runs over $(0,\infty)$.
When the type needs not be specified, we simply write $f\sim_{\{\M\}}^u\widehat{f}$ in $S$.
Note that, taking $p=0$ in~\eqref{desarasintunifo}, we deduce that every function in $\widetilde{\mathcal{A}}^u_{\{\M\}}(S)$ is a bounded function.

Finally, we define for every $A>0$ the class $\mathcal{A}_{\M,A}(S)$ consisting of the holomorphic functions $f$ in $S$ such that
$$
\left\|f\right\|_{\M,A}:=\sup_{z\in S,p\in\N_{0}}\frac{|f^{(p)}(z)|}{A^{p}M_{p}}<\infty.
$$
($\mathcal{A}_{\M,A}(S),\left\|\,\cdot\, \right\| _{\M,A}$) is a Banach space, and $\mathcal{A}_{\{\M\}}(S):=\cup_{A>0}\mathcal{A}_{\M,A}(S)$ is called a \emph{Carleman-Roumieu ultraholomorphic class} in the sector $S$, whose natural inductive topology makes it an $(LB)$ space.\par
We warn the reader that these notations, while the same as in the paper~\cite{JimenezSanzSchindlSurjectDC}, do not agree with the ones used in~\cite{SanzFlatProxOrder,JimenezSanzSchindlInjectSurject}, where
$\widetilde{\mathcal{A}}^u_{\{\M\}}(S)$ was denoted by
$\widetilde{\mathcal{A}}^u_{\M}(S)$, $\mathcal{A}_{\M,A}(S)$ by $\mathcal{A}_{\M/\bL_1,A}(S)$, and $\mathcal{A}_{\{\M\}}(S)$ by $\mathcal{A}_{\M/\bL_1}(S)$.

If $\M$ is (lc), the spaces $\mathcal{A}_{\{\M\}}(S)$ and $\widetilde{\mathcal{A}}^u_{\{\M\}}(S)$
are algebras, and if $\M$ is (dc) they are stable under taking derivatives.
Moreover, if $\M\approx\bL$ the corresponding classes coincide.

Since the derivatives of $f\in\mathcal{A}_{\M,A}(S)$ are Lipschitz, for every $p\in\N_{0}$ one may define
\begin{equation}\label{eq.deriv.at.0.def}
f^{(p)}(0):=\lim_{z\in S,z\to0 }f^{(p)}(z)\in\C.
\end{equation}

As a consequence of Taylor's formula and Cauchy's integral formula for the derivatives, there is a close relation between Carleman-Roumieu ultraholomorphic classes and the concept of asymptotic expansion (the proof may be easily adapted from~\cite{balserutx}).

\begin{prop}\label{propcotaderidesaasin}
Let $\M$ be a sequence and $S$ be a sector. Then,
\begin{enumerate}[(i)]
\item If $f\in\mathcal{A}_{\hM,A}(S)$ then $f$ admits $\widehat{f}:=\sum_{p\in\N_0}\frac{1}{p!}f^{(p)}(0)z^p$ as its uniform $\M$-asymptotic expansion in $S$ of type $1/A$, where $(f^{(p)}(0))_{p\in\N_0}$ is given by \eqref{eq.deriv.at.0.def}. Moreover, $\|f\|_{\M,A,\overset{\sim}{u}}\le \|f\|_{\hM,A}$, and so the identity map $\mathcal{A}_{\hM,A}(S)\hookrightarrow \widetilde{\mathcal{A}}^u_{\M,A}(S)$ is continuous. Consequently, we also have that
$\mathcal{A}_{\{\hM\}}(S)\subseteq \widetilde{\mathcal{A}}^u_{\{\M\}}(S)$
and $\mathcal{A}_{\{\hM\}}(S)\hookrightarrow \widetilde{\mathcal{A}}^u_{\{\M\}}(S)$ is continuous.

\item If $S$ is unbounded and $T$ is a proper subsector of $S$, then there exists a constant $c=c(T,S)>0$ such that the restriction to $T$, $f|_T$, of functions $f$ defined on $S$ and admitting a uniform $\M$-asymptotic expansion in $S$ of type $1/A>0$, belongs to $\mathcal{A}_{\hM,cA}(T)$, and $\|f|_T\|_{\hM,cA}\le \|f\|_{\M,A,\overset{\sim}{u}}$. So, the restriction map from $\widetilde{\mathcal{A}}^u_{\M,A}(S)$ to $\mathcal{A}_{\hM,cA}(T)$ is continuous, and it is also continuous from $\widetilde{\mathcal{A}}^u_{\{\M\}}(S)$ to $\mathcal{A}_{\{\hM\}}(T)$.
\end{enumerate}
\end{prop}

One may accordingly define classes of formal power series
\begin{equation}\label{eq.defBanachFormalPowerSeries}
\C[[z]]_{\M,A}=\Big\{\widehat{f}=\sum_{p=0}^\infty a_pz^p\in\C[[z]]:\, \left|\,\widehat{f} \,\right|_{\M,A}:=\sup_{p\in\N_{0}}\displaystyle \frac{|a_{p}|}{A^{p}M_{p}}<\infty\Big\}.
\end{equation}%
$(\C[[z]]_{\M,A},\left| \,\cdot\,  \right|_{\M,A})$ is a Banach space and we put $\C[[z]]_{\{\M\}}:=\cup_{A>0}\C[[z]]_{\M,A}$, again an $(LB)$ space.

It is natural to consider the \emph{asymptotic Borel map}  $\widetilde{\mathcal{B}}$
sending a function $f\in\widetilde{\mathcal{A}}^u_{\M,A}(S)$
into its $\M$-asymptotic expansion $\widehat{f}\in\C[[z]]_{\M,A}$.
By Proposition~\ref{propcotaderidesaasin}.(i)
the asymptotic Borel map may be defined from $\widetilde{\mathcal{A}}^u_{\{\M\}}(S)$ or $\mathcal{A}_{\{\hM\}}(S)$ into  $\C[[z]]_{\{\M\}}$, and from $\mathcal{A}_{\hM,A}(S)$ into $\C[[z]]_{\M,A}$.

If $\M$ is (lc), $\widetilde{\mathcal{B}}$ is a homomorphism of algebras; if $\M$ is also (dc), differentiation commutes with $\widetilde{\mathcal{B}}$. Moreover, it is continuous when considered between the corresponding Banach or $(LB)$ spaces previously introduced. Finally, note that if $\M\approx\bL$, then $\C[[z]]_{\{\M\}}=\C[[z]]_{\{\bL\}}$, and the corresponding Borel maps are in all cases identical.

Since the problem under study is invariant under rotation, we will focus on the surjectivity of the Borel map in unbounded sectors $S_{\gamma}$. So, we define
\begin{align*}
S_{\{\hM\}}:=&\{\gamma>0; \quad \widetilde{\mathcal{B}}:\mathcal{A}_{\{\hM\}}(S_\gamma)\longrightarrow \C[[z]]_{\{\M\}} \text{ is surjective}\} ,\\
\widetilde{S}^u_{\{\M\}}:=&\{\gamma>0; \quad\widetilde{\mathcal{B}}:\widetilde{\mathcal{A}}^u_{\{\M\}}(S_\gamma)\longrightarrow \C[[z]]_{\{\M\}} \text{ is surjective}\}.
\end{align*}
We again note that these intervals were respectively denoted by $S_{\M}$ and $\widetilde{S}^u_{\M}$
in~\cite{JimenezSanzSchindlInjectSurject}.\par

It is clear that $S_{\{\hM\}}$ and $\widetilde{S}^u_{\{\M\}}$
are either empty or left-open intervals  having $0$ as endpoint, called \emph{surjectivity intervals}.
Using~Proposition~\ref{propcotaderidesaasin},
we easily see that
\begin{align}
(\widetilde{S}^u_{\{\M\}})^{\circ}\subseteq S_{\{\hM\}}\subseteq \widetilde{S}^u_{\{\M\}},
\label{equaContentionSurjectIntervals}
\end{align}
where $I^{\circ}$ stands for the interior of the interval $I$.

\section{Optimal flat functions and surjectivity of the Borel map}\label{sectFlatFunctions}

The following result appeared, in a slightly different form, in~\cite[Lemma~4.5]{JimenezSanzSchindlInjectSurject}.

\begin{lemma}\label{lemma.SurjectImpliesgammaMpositive}
Let $\M$ be a weight sequence. If $\widetilde{S}^u_{\{\M\}}\neq\emptyset$, then $\M$ satisfies (snq) or, equivalently, $\ga(\M) >0$.
\end{lemma}

Subsequently, a converse, more precise statement appeared in~\cite[Th.~3.7]{JimenezSanzSchindlSurjectDC} under the additional hypothesis of condition (dc).

\begin{theo}\label{teorSurject.dc}
Let $\hM$ be a regular sequence such that $\gamma(\M)>0$. Then,
$$(0,\gamma(\M))\subseteq S_{\{\hM\}}\subseteq \widetilde{S}^u_{\{\M\}}\subseteq
(0,\gamma(\M)].
$$
In particular, if $\gamma(\M)= \infty$, we have that $S_{\{\hM\}}= \widetilde{S}^u_{\{\M\}}= (0,\infty)$.
\end{theo}

So, the surjectivity of the Borel map for regular sequences is governed by the value of the index $\ga(\M)$.

\subsection{Construction of optimal flat functions}

Our aim is to relate the surjectivity of the Borel map in a sector to the existence of optimal flat functions in it, which we now define and construct in this subsection.

\begin{defi}\label{optimalflatdef}
Let $\M$ be a weight sequence, $S$ an unbounded sector bisected by direction $d=0$, i.e., by the positive real line $(0,+\infty)\subset\mathcal{R}$. A holomorphic function $G\colon S\to\C$ is called an \emph{optimal $\{\M\}$-flat function} in $S$ if:
	\begin{itemize}
		\item[$(i)$] There exist $K_1,K_2>0$ such that for all $x>0$,
		\begin{equation}\label{optimalflatleft}
			K_1h_{\M}(K_2x)\le G(x).
		\end{equation}
		
		\item[$(ii)$] There exist $K_3,K_4>0$ such that for all $z\in S$, one has
		\begin{equation}\label{optimalflatright}
			|G(z)|\le K_3h_{\M}(K_4|z|).
		\end{equation}
	\end{itemize}
\end{defi}

Besides the symmetry imposed by condition $(i)$ (observe that $G(x)>0$ for $x>0$, and so $G(\overline{z})=\overline{G(z)}$, $z\in S$),
we note that the estimates in~\eqref{optimalflatright} amount to the fact that
$$|G(z)|\le K_3K_4^pM_p|z|^p,\qquad p\in\N_0,\ z\in S,
$$
what exactly means that $G\in\widetilde{\mathcal{A}}^u_{\{\M\}}(S)$ and is \emph{$\{\M\}$-flat}, i.e., its uniform $\{\M\}$-asymptotic expansion is given by the null series. The inequality imposed in~\eqref{optimalflatleft} makes the function optimal in a sense, as its rate of decrease on the positive real axis when $t$ tends to 0 is accurately specified by the function $h_{\M}$.
Note that, in previous instances where such optimal flat functions appear~\cite{Thilliez03,LastraMalekSanzContinuousRightLaplace,JimenezLastraSanzRapidGrowth}, the estimates from below in~\eqref{optimalflatleft} are imposed and/or obtained in the whole sector $S$, and not just on its bisecting direction. We think the present definition is more convenient, since it is easier to check for concrete functions, and for our purposes it provides all the necessary information in order to work with such functions.


As a first step for the construction of such flat functions, we need to estimate the harmonic extension $P_\sigma$ in terms of the majorant $\ka_\sigma$. The right-hand side estimate in the next result is a slight refinement of the one in~\cite[Lemma\ 3.3]{BonetMeiseTaylor92}, which was not precise enough for our purposes. We include the whole proof for the sake of completeness.

\begin{prop}
Let $\sigma:[0,\infty)\to [0,\infty)$ be a nondecreasing nonquasianalytic function. Then, we have
	\begin{equation}\label{kappaPbound}
		\frac{1}{\pi}\ka_\sigma(y)\leq P_\sigma(iy)\leq \ka_\sigma(y),\qquad y\geq0.
	\end{equation}
\end{prop}
\begin{proof}
	If $y=0$ all the values are equal to $\sigma(0)$ and so the inequalities hold true. Now, for $y>0$ we have
\begin{equation*}
P_\sigma(iy)=\frac{y}{\pi}\int_{-\infty}^{\infty}\frac{\sigma(|t|)}{t^2+y^2}dt=
\frac{2y}{\pi}\int_{0}^{\infty}\frac{\sigma(t)}{t^2+y^2}dt=\frac{2}{\pi}\int_{0}^{\infty}\frac{\sigma(ys)}{s^2+1}ds\geq \frac{2}{\pi}\int_{1}^{\infty}\frac{\sigma(ys)}{s^2+1}ds.
\end{equation*}
%
Since $s^2+1\leq 2s^2$ for $s\geq 1$, we deduce that
	\begin{equation*}
		P_\sigma(iy)\geq 
\frac{1}{\pi}\int_{1}^{\infty}\frac{\sigma(ys)}{s^2}ds=\frac{1}{\pi}\ka_\sigma(y).
	\end{equation*}
In order to prove the right inequality, we start by splitting the integral into two parts:
	\begin{equation}\label{splitintegralPsigma}
P_\sigma(iy)=\frac{2}{\pi}\int_{0}^{\infty}\frac{\sigma(ys)}{s^2+1}ds= \frac{2}{\pi}\left(\int_{0}^{1}\frac{\sigma(ys)}{s^2+1}ds + \int_{1}^{\infty}\frac{\sigma(ys)}{s^2+1}ds\right).
	\end{equation}
As $\sigma$ is nondecreasing, we may write
	\begin{equation}\label{splitintegralPsigma1}
\int_{0}^{1}\frac{\sigma(ys)}{s^2+1}ds\leq\sigma(y)\int_{0}^{1}\frac{1}{s^2+1}ds=\frac{\pi}{4}\sigma(y)
	\end{equation}
and
	\begin{equation}\label{splitintegralPsigma2}
\int_{1}^{\infty}\frac{\sigma(ys)}{s^2+1}ds=\ka_\sigma(y)-\int_{1}^{\infty}\left( \frac{1}{s^2}-\frac{1}{s^2+1}\right)\sigma(ys)ds\leq  \ka_\sigma(y)- \sigma(y)\left(1-\frac{\pi}{4}\right).
	\end{equation}
From \eqref{splitintegralPsigma}, \eqref{splitintegralPsigma1}, \eqref{splitintegralPsigma2} and \eqref{sigmakappabound} we deduce that
	\begin{equation*}
		P_\sigma(iy)\leq \frac{2}{\pi}\left(\frac{\pi}{2}\sigma(y) +\ka_\sigma(y)- \sigma(y)\right)\leq \ka_\sigma(y).
	\end{equation*}
\end{proof}

The key condition for weight sequences that will allow us to construct optimal flat functions appeared in a work of M.~Langenbruch~\cite{Langenbruch94}.

\begin{defi} Let $\M$ be a weight sequence such that $\widecheck{\M}$ satisfies (nq), so that $P_{\omega_{\M}}$ is well-defined. We say that the sequence satisfies the \emph{Langenbruch's condition} if there exists a constant $C>0$ such that for all $y\geq0$ we have
\begin{equation}\label{eq.Langenbruch}
		P_{\omega_{\M}}(iy)\leq\omega_\M(Cy)+C.
\end{equation}
\end{defi}

We can characterize the previous condition in terms of the index $\gamma(\M)$. This connection has very recently appeared for the first time in a work of D. N. Nenning, A. Rainer and the fourth author~\cite{NenningRainerSchindl_preprint}. Although the additional hypothesis of (dc) appears in their (indirect) arguments, it can be removed as long as the sequence satisfies (snq), as we now show. Observe that, by Lemma~\ref{lemma.SurjectImpliesgammaMpositive}, the condition (snq) (equivalently, $\gamma(\M)>0$) is necessary for surjectivity, so it is not a restriction for our aim.

%
%
%
%
%
%
%
%
%
%

\begin{prop}\label{langenbruch} Let $\M$ be a weight sequence. The following are equivalent:
	\begin{itemize}
		\item[(i)] $\gamma(\M)>0$, $\widecheck{\M}$ satisfies (nq) and $\M$ satisfies Langenbruch's condition.
		\item[(ii)] $\gamma(\M)>1$.
	\end{itemize}
\end{prop}
\begin{proof}
	 First, from \eqref{assofuncintegral} we deduce that for all $r\geq 0$ and $B\geq 0$,
	\begin{equation}\label{omegaBrbound}
		\omega_{\M}(e^Br)=\int_0^{e^B r}\frac{\nu_{\m}(u)}{u}du=\omega_{\M}(r)+\int_r^{e^B r}\frac{\nu_{\m}(u)}{u}du\geq \omega_{\M}(r)+B\nu_{\m}(r).
	\end{equation}
The last inequality is a consequence of the monotonicity of $\nu_{\m}$.

(i)$\Rightarrow$(ii) By taking into account \eqref{omegapbound} and \eqref{kappaPbound}, we deduce
\begin{equation*}
	\omega_{\M}(y)+\kappa_{\nu_{\m}}(y)\leq  P_{\omega_{\M}}(iy)+\pi P_{\nu_{\m}}(iy)= P_{\omega_{\M}+\pi\nu_{\m}}(iy),	\qquad y\geq 0.
	\end{equation*}
	Thanks to \eqref{omegaBrbound} and the monotonicity of the harmonic extension with respect to the argument function we get from above
	\begin{equation*}
		\omega_{\M}(y)+\kappa_{\nu_{\m}}(y)\leq P_{\omega_{\M}(e^\pi\cdot)}=P_{\omega_{\M}}(ie^\pi y)\leq\omega_{\M}(Ce^\pi y)+C, \qquad y\geq 0.
	\end{equation*}
Next, by using the integral expression \eqref{assofuncintegral} and the monotonicity of $\nu_{\m}$ we have that
	\begin{equation*}
		\kappa_{\nu_{\m}}(y)\leq \omega_{\M}(Ce^\pi y)-\omega_{\M}(y)+C=\int_y^{Ce^\pi y}\frac{\nu_{\m}(u)}{u}du+C	 \leq \ln(Ce^\pi)\nu_{\m}(Ce^\pi y)+C, \qquad y\geq 0.
	\end{equation*}
	Finally, by Lemma~\ref{functionNUproperties}, we deduce that
	\begin{equation*}
		\kappa_{\nu_{\m}}(y)\leq D\nu_{\m}(y)+D, \qquad y\geq 0,
	\end{equation*}
for suitable $ D>0$. This is condition $(\omega_{\text{snq}})$ for $\nu_{\m}$ and, by Lemma~\ref{functionNUproperties}, we may conclude that $\gamma(\M)>1$.\\

(ii)$\Rightarrow$(i) Condition $\gamma(\M)>1$ implies that $\gamma(\M)>0$, and amounts to condition $(\gamma_1)$ for $\bm$ (see~\eqref{equa.indice.gammaM.gamma_r}), so that $\widecheck{\M}$ clearly satisfies (nq). By Lemma~\ref{functionNUproperties}, the condition $\gamma(\M)>1$ is equivalent to the existence of a constant $C>0$ such that
	\begin{equation}\label{kappanu-komatsu}
		\kappa_{\nu_{\m}}(y)\leq C\nu_{\m}(y)+C, \qquad y\geq 0.
	\end{equation}
	Then, from \eqref{kappaPbound}, \eqref{kappaomeganurelation} and the above inequality we deduce that
	\begin{equation*}
		P_{\omega_{\M}}(iy) \leq \kappa_{\omega_{\M}}(y)=\omega_{\M}(y)+\kappa_{\nu_{\m}}(y) \leq\omega_{\M}(y)+C\nu_{\m}(y)+C, \qquad y\geq 0.
	\end{equation*}
	By \eqref{omegaBrbound}, we have from above that
	\begin{equation*}
		P_{\omega_{\M}}(iy) \leq \omega_{\M}(e^Cy)+C, \qquad y\geq 0,
	\end{equation*}
which completes the proof.
\end{proof}

\begin{rema}
The condition $\gamma(\M)>1$ is the same as $\gamma(\widecheck{\M})>0$, or equivalently, (snq) for $\widecheck{\M}$ (even if $\widecheck{\M}$ might not satisfy (lc), we can apply~\cite[Coro.\ 3.13]{JimenezSanzSchindlIndex} to obtain this equivalence). So, for a weight sequence $\M$ satisfying (snq), Langenbruch's condition allows to pass from (nq) to (snq) for $\widecheck{\M}$.

Observe also that, by~\cite[Lemma\ 3.20]{JimenezSanzSchindlIndex}, the condition (nq) for $\widecheck{\M}$ implies that the index $\omega(\widecheck{\M})$, introduced in~\cite{SanzFlatProxOrder} and studied in detail in~\cite{JimenezSanzSchindlIndex}, is nonnegative, and so $\omega(\M)=\omega(\widecheck{\M})+1\ge 1$. As one only knows that $\ga(\M)\le\omega(\M)$ in general, and these indices can perfectly be different, one may better understand the effect of Langenbruch's condition.
\end{rema}

\begin{rema}
On the one hand, as said before, for a weight sequence $\M$ the condition $\ga(\M)>1$ amounts to the condition $(\ga_1)$ for $\bm$, and it is well-known (see~\cite[Prop.\ 4.4]{komatsu}) that then $\o_{\M}$ satisfies $(\omega_{\text{snq}})$.  As it can be deduced from~\cite[Prop.\ 1.7]{MeiseTaylor88}, this last fact is, in its turn, equivalent to the existence of a constant $C>0$ such that
$$
P_{\o_{\M}}(iy)\le C\o_{\M}(y)+C,\quad y\ge 0.
$$
On the other hand, in~\cite[Prop.\ 3.6]{komatsu} the condition (mg) for a weight sequence $\M$ is shown to be equivalent to the fact that $2\o_{\M}(y)\le \o_{\M}(Dy)+D$ for all $y\ge 0$ and suitable $D>0$. Gathering these estimates, we conclude that if $\M$ is strongly regular then $\ga(\M)>1$ if, and only if, $\M$ satisfies Langenbruch's condition. This was basically the reasoning that allowed V. Thilliez to obtain optimal $\{\M\}$-flat functions, in the very same way as we are doing in the next result, but dropping now the moderate growth condition by means of Proposition~\ref{langenbruch}.
\end{rema}

Thanks to the previous result, we will construct optimal $\{\M\}$-flat functions in the right half plane as long as $\gamma(\M)>1$.

\begin{prop}\label{optimalff-S1}	Let $\M$ be a weight sequence. If $\ga(\M)>1$, then the function $$
G(z)=\exp(-P_{\omega_\M}(i/z)-iQ_{\omega_\M}(i/z))
$$
is an optimal $\{\M\}$-flat function in the halfplane $S_1$, where $Q_{\omega_\M}$ is the harmonic conjugate of $P_{\omega_\M}$ in the upper half plane.
\end{prop}
\begin{proof}
It is clear that the function $G$ is holomorphic in $S_1$. On the one hand, by taking into account \eqref{omegapbound}, for $z\in S_1$ we have that
	$$|G(z)|=\exp(-P_{\omega_\M}(i/z))\leq \exp(-\omega_\M(1/|z|))=h_\M(|z|).$$
On the other hand, the condition $\ga(\M)>1$ implies, by Proposition \ref{langenbruch}, that there exists $C>0$ such that $P_{\omega_{\M}}(ix)\leq\omega_\M(Cx)+C$ for every $x>0$. Since one can easily check that $Q_{\omega_\M}(i/x)=0$, we have that
$$G(x)=\exp(-P_{\omega_\M}(i/x))\geq\exp(-\omega_\M(C/x)-C)= \exp(-C)h_\M(x/C),$$
as desired.
\end{proof}

By a ramification of the variable we can extend this method to an arbitrary weight sequence with $\ga(\M)>0$ and any sector whose opening is less than $\pi\ga(\M)$.

\begin{prop}\label{prop.ExistOptFlatFunct}
Let $\M$ be a weight sequence with $\ga(\M)>0$. Then, for any $0<\gamma<\ga(\M)$ there exist an optimal $\{\M\}$-flat function in $S_\ga$.
\end{prop}
\begin{proof}
	Let $s>0$ be such that $\ga<1/s<\ga(\M)$. Then, by \cite[Th.\ 3.10, Prop.\ 3.6]{JimenezSanzSchindlIndex} we have that $\ga({\M}^s)=s\ga(\M)>1$, where ${\M}^s:=(M_p^s)_{p\in\N_0}$ is again a weight sequence. We apply the last result to the sequence ${\M}^s$, so there exist an optimal $\{\M^s\}$-flat function $G$ in $S_1$. It is important to note that the bounds for $G$ appearing in Definition \ref{optimalflatdef} will be in terms of $h_{{\M}^{s}}$, instead of $h_{\M}$. Moreover, the following relation between the functions $\omega_{\M^s}$ and $\omega_{\M}$ is straightforward:
	\begin{equation}\label{omega-Ms}
		\omega_{\M}(t^{1/s})=\frac{1}{s}\omega_{\M^s}(t),\qquad t\geq 0.
	\end{equation}
	Now, let us prove that the function $F(z)=(G(z^s))^{1/s}$, $z\in S_\ga$, is an optimal $\{\M\}$-flat function in $S_\ga$.
	From the fact that $G$ is an optimal $\{\M^s\}$-flat function, \eqref{functionhequ2} and \eqref{omega-Ms}, we get
		\begin{align*}
				F(x)=(G(x^s))^{1/s}&\geq K_1^{1/s}\exp(-s^{-1}\omega_{\M^s}(1/(K_2 x^s)))\\
				&\geq K_1^{1/s}\exp(-\omega_{\M}(1/(K_2^{1/s} x)))=K_1^{1/s}h_{\M}(K_2^{1/s} x), \qquad x>0,
		\end{align*}
	for suitable constants $K_1,K_2>0$. Moreover, we have that
		\begin{align*}
				|F(z)|&\leq K_3^{1/s}\exp(-s^{-1}\omega_{\M^s}(1/(K_4 |z|^s)))\\
				&\leq K_3^{1/s}\exp(-\omega_{\M}(1/(K_4^{1/s} |z|)))=K_3^{1/s}h_{\M}(K_4^{1/s} |z|),\qquad z\in S_\ga,
		\end{align*}
	for suitable constants $K_3,K_4>0$, and we are done.
\end{proof}

\subsection{Surjectivity of the Borel map for regular sequences}

We will describe next how, by means of an optimal flat function, one can obtain extension operators, right inverses for the Borel map, for ultraholomorphic classes defined by regular sequences.

If $G$ is an optimal $\{\M\}$-flat function in $\widetilde{\mathcal{A}}^u_{\{\M\}}(S)$, we define the kernel function $e\colon S\to\C$ given by
\begin{equation*}
e(z):=G\left(\frac{1}{z}\right),\quad z\in S.
\end{equation*}
It is obvious that $e(x)>0$ for all $x>0$, and there exist $K_1,K_2,K_3,K_4>0$ such that
\begin{equation}\label{eq.Bounds_e_sector}
K_1h_{\M}\left(\frac{K_2}{x}\right)\le e(x), \quad x>0, \hspace{1cm}\text{and}\hspace{1cm} |e(z)|\le K_3h_{\M}\left(\frac{K_4}{|z|}\right),\quad z\in S.
\end{equation}
For every $p\in\mathbb{N}_0$ we define the \emph{$p$-th moment} of the function $e(z)$, given by
\begin{equation*}
\mo(p):=\int_0^\infty t^{p}e(t)\,dt.
\end{equation*}
Note that the positive value $\mo(0)$ need not be equal to 1.

The following result is crucial for our aim.

\begin{prop}\label{prop.boundsMomentsSigmaEntre1y2}
Suppose $\hM$ is a regular sequence and $G$ is an optimal $\{\M\}$-flat function in $\widetilde{\mathcal{A}}^u_{\{\M\}}(S)$. Consider the sequence of moments $\mo:=(\mo(p))_{p\in\N_0}$ associated with the kernel function $e(z)=G(1/z)$. Then, there exist $B_1,B_2>0$ such that
\begin{equation}\label{eq.boundsMmomentsdc}
\mo(0)B_1^{p}M_p\le \mo(p)\le \mo(0)B_2^{p}M_p,\quad p\in\mathbb{N}_0.
\end{equation}
In other words, $\M$ and $\mo$ are equivalent.
\end{prop}

\begin{proof1}
We only need to reason for $p\in\N$.
On the one hand, because of the right-hand inequalities in~\eqref{eq.Bounds_e_sector} and Lemma~\ref{functionhproperties}.(ii), for every $p\in\N$ and $s>0$ we may write
\begin{align*}
\mo(p)&=\int_0^s t^p e(t)\,dt
+\int_s^\infty \frac{1}{t^2}t^{p+2}e(t)\,dt\\
&\le K_3 \int_0^s t^p\,dt + K_3\sup_{t>0}t^{p+2}h_\M\left(\frac{K_4}{t}\right)\int_s^\infty \frac{1}{t^2}\,dt \\
&= K_3 \frac{s^{p+1}}{p+1}
+K_3\frac{1}{s}K_4^{p+2}M_{p+2}
\le K_3 \left(\frac{s^{p+1}}{p+1}
+\frac{(K_4D)^{p+2}M_{p}}{s}\right).
\end{align*}
Note that in the last equality we have used \eqref{eq.MpfromomegaM}, and then we have applied (dc) with a suitable constant $D>0$.
Since $s>0$ was arbitrary, we finally get
\begin{equation*}
\mo(p)\le \inf_{s>0}K_3 \left(\frac{s^{p+1}}{p+1}
+\frac{(K_4D)^{p+2}M_{p}}{s}\right)= K_3\frac{p+2}{p+1}(K_4D)^{p+1}(M_{p})^{(p+1)/(p+2)}
\le \mo(0)B_2^pM_p
\end{equation*}
for a suitably enlarged constant $B_2>0$ (observe that $p\ge 1$ and that, eventually, $M_p\ge 1$).

On the other hand, by the left-hand inequalities in~\eqref{eq.Bounds_e_sector} and Lemma~\ref{functionhproperties}.(i), for every $p\in\N$ and $s>0$ we may estimate
$$
\mo(p)\ge \int_0^s t^p e(t)\,dt \ge K_1 \int_0^s t^p h_{\M}\left(\frac{K_2}{t}\right)\,dt\ge
K_1 h_{\M}\left(\frac{K_2}{s}\right)\frac{s^{p+1}}{p+1}.
$$
Then, again by~\eqref{eq.MpfromomegaM}, we deduce that
$$
\mo(p)\ge \frac{K_1}{p+1} \sup_{s>0}h_{\M}\left(\frac{K_2}{s}\right)s^{p+1}=
\frac{K_1}{p+1} K_2^{p+1}M_{p+1}\ge
\mo(0)B_1^p M_p
$$
for a suitable constant $B_1>0$ (note that $\M$ is eventually nondecreasing).
\end{proof1}

We can already state the following main result. The forthcoming implication $(ii)\Rightarrow(v)$ for strongly regular sequences $\M$ was first obtained by V. Thilliez~\cite[Th.~3.2.1]{Thilliez03}, and the proof heavily rested on the moderate growth condition,
both for the  construction~\cite[Th.~2.3.1]{Thilliez03} of optimal $\{\M\}$-flat functions in sectors $S_{\ga}$ for every $\ga>0$ such that $\ga<\ga(\M)$, and for the subsequent use
of Whitney extension results in the ultradifferentiable setting.
In~\cite{LastraMalekSanzContinuousRightLaplace} the implication $(ii)\Rightarrow(iii)$ was proved again for strongly regular sequences, but with a completely different technique, and it is this approach which allows here for the weakening of condition (mg) into (dc).

\begin{theo}\label{tpral}
Let $\hM$ be a regular sequence (that is, $\M$ is a weight sequence and satisfies (dc))
with $\ga(\M)>0$, and let $\ga>0$ be given. Each of the following statements implies the next one:
\begin{itemize}
\item[(i)] $\ga<\ga(\M)$.
\item[(ii)] There exists an optimal $\{\M\}$-flat function in $\widetilde{\mathcal{A}}^u_{\{\M\}}(S_{\ga})$.
\item[(iii)] There exists $c>0$ such that for every $A>0$ there exists a linear continuous map $T_{\M,A}\colon\C[[z]]_{\M,A}\to \widetilde{\mathcal{A}}^u_{\M,cA}(S_{\ga})$ such that $\widetilde{\mathcal{B}}\circ T_{\M,A}$ is the identity map in $\C[[z]]_{\M,A}$ (i.e., $T_{\M,A}$ is an extension operator,  right inverse for $\widetilde{\mathcal{B}}$).
\item[(iv)] The Borel map $\widetilde{\mathcal{B}}\colon \widetilde{\mathcal{A}}^u_{\{\M\}}(S_{\ga})\to\C[[z]]_{\{\M\}}$ is surjective. In other words, $(0,\ga]\subset\widetilde{S}_{\{\M\}}^u$.
\item[(v)] $(0,\ga)\subset S_{\{\hM\}}$.
\item[(vi)] $\ga\le\ga(\M)$.
\end{itemize}
\end{theo}

\begin{proof1}
$(i)\Rightarrow (ii)$ See Proposition~\ref{prop.ExistOptFlatFunct}, valid for any weight sequence $\M$.

$(ii)\Rightarrow (iii)$ Let $A>0$ and $\widehat{f}=\sum_{p=0}^\infty a_pz^p\in\C[[z]]_{\M,A}$ be given. Let $(\mo(p))_{p\in\N_0}$ be the sequence of moments associated to the function $e(z)=G(1/z)$, where $G$ is an optimal $\{\M\}$-flat function in $\widetilde{\mathcal{A}}^u_{\{\M\}}(S_{\ga})$.
By the definition of the norm in $\C[[z]]_{\M,A}$ (see~\eqref{eq.defBanachFormalPowerSeries}), we have
\begin{equation*}
|a_p|\le |\widehat{f}|_{\M,A}A^{p}M_p,\quad p\in\N_0.
\end{equation*}
From the left-hand inequalities in~\eqref{eq.boundsMmomentsdc}, we deduce that
\begin{equation}\label{eq.boundsCoeffBorelTransf}
\left|\frac{a_p}{\mo(p)}\right|\le \frac{|\widehat{f}|_{\M,A}}{\mo(0)}\left(\frac{A}{B_1}\right)^{p},\quad p\in\N_0.
\end{equation}
Hence, the formal Borel-like transform of $\widehat{f}$,
$$
\widehat{g}=\sum_{p=0}^\infty\frac{a_p}{\mo(p)}z^p,
$$
is convergent in the disc $D(0,R)$ for $R=B_1/A>0$, and it defines a holomorphic function $g$ there. Choose $R_0:=B_1/(2A)<R$, and define
\begin{equation*}
T_{\M,A}(\widehat{f}\,)(z):=\frac{1}{z}\int_{0}^{R_0}e\left(\frac{u}{z}\right)g(u)\,du,\qquad z\in S_{\ga},
\end{equation*}
which is a truncated Laplace-like transform of $g$ with kernel $e$.
By virtue of Leibniz's theorem for parametric integrals and the properties of $e$, we deduce that this function, denoted by $f$ for the sake of brevity, is holomorphic in $S_{\ga}$. We will prove that
$f\sim^u_{\{\bM\}}\hat{f}$ uniformly in $S_{\ga}$, and that
the map $\widehat{f}\mapsto f$, which is obviously linear, is also continuous from $\C[[z]]_{\M,A}$ into $\widetilde{\mathcal{A}}^u_{\M,cA}(S_{\ga})$ for suitable $c>0$ independent from $A$.

Let $p\in\N_0$ and $z\in S_{\ga}$. We have
\begin{align*}
f(z)-\sum_{n=0}^{p-1}a_nz^n &= f(z)-\sum_{n=0}^{p-1}\frac{a_n}{\mo(n)}\mo(n)z^n\\
&= \frac{1}{z}\int_{0}^{R_0}e \left(\frac{u}{z}\right) \sum_{n=0}^{\infty}\frac{a_{n}}{\mo(n)}u^n\,du -\sum_{n=0}^{p-1}\frac{a_n}{\mo(n)}\int_{0}^{\infty}v^{n}e(v)\,dv\, z^n.
\end{align*}
After a change of variable $u=zv$ in the last integral, one may use Cauchy's residue theorem and the right-hand estimates in~(\ref{eq.Bounds_e_sector}) in order to rotate the path of integration and obtain
$$
z^n\int_{0}^{\infty}v^{n}e(v)dv= \frac{1}{z}\int_{0}^{\infty}u^{n}e\left(\frac{u}{z}\right)\,du.
$$
So, we can write the preceding difference as
\begin{equation*}
\frac{1}{z}\left(\int_{0}^{R_0}e\left(\frac{u}{z}\right) \sum_{n=p}^{\infty}\frac{a_{n}}{\mo(n)}u^n\,du -\int_{R_0}^{\infty}e\left(\frac{u}{z}\right) \sum_{n=0}^{p-1}\frac{a_n}{\mo(n)}u^{n}\,du\right).
\end{equation*}
Then, we have
\begin{equation}\label{eq.RemainderTwoSums}
\left|f(z)-\sum_{n=0}^{p-1}a_nz^n\right|\le \frac{1}{|z|}(f_{1}(z)+f_2(z)),
\end{equation}
where
$$
f_{1}(z)=\left|\int_{0}^{R_0}e \left(\frac{u}{z}\right) \sum_{n=p}^{\infty}\frac{a_{n}}{\mo(n)}u^n \,du\right|,\quad
f_{2}(z)=\left|\int_{R_0}^{\infty}e \left(\frac{u}{z}\right) \sum_{n=0}^{p-1}\frac{a_n}{\mo(n)}u^n\,du\right|.$$
We now estimate $f_1(z)$ and $f_2(z)$.
Observe that for every $u\in(0,R_0]$ we have $0<Au/B_1\le 1/2$. So, from~\eqref{eq.boundsCoeffBorelTransf} we get
\begin{equation*}
\sum_{n=p}^{\infty}\frac{|a_{n}|}{\mo(n)}u^n\le \frac{|\widehat{f}|_{\M,A}}{\mo(0)} \sum_{n=p}^{\infty}\left(\frac{Au}{B_1}\right)^n\le
\frac{2|\widehat{f}|_{\M,A}}{\mo(0)}\left(\frac{A}{B_1}\right)^p u^p.
\end{equation*}
Hence,
\begin{equation}\label{eq.Estimates_f1}
f_{1}(z)\le \frac{2|\widehat{f}|_{\M,A}}{\mo(0)} \left(\frac{A}{B_1}\right)^p \int_{0}^{R_0}\left|e\left(\frac{u}{z}\right)\right| u^p\,du.
\end{equation}
Regarding $f_{2}(z)$, for $u\ge R_0$ and $0\le n\le p-1$ we have $(u/R_0)^n\le (u/R_0)^p$, so $u^n\le R_0^nu^p/R_0^p$. Again by~(\ref{eq.boundsCoeffBorelTransf}), and taking into account the value of $R_0$, we may write
$$\sum_{n=0}^{p-1}\frac{|a_n|}{\mo(n)}u^n\le \frac{|\widehat{f}|_{\M,A}}{\mo(0)}\frac{u^p}{R_0^p} \sum_{n=0}^{p-1}\left(\frac{AR_0}{B_1}\right)^n\le \frac{|\widehat{f}|_{\M,A}}{\mo(0)}\left(\frac{2A}{B_1}\right)^p u^p.$$
Then, we get
\begin{equation}\label{eq.Estimates_f2}
f_2(z)\le \frac{|\widehat{f}|_{\M,A}}{\mo(0)} \left(\frac{2A}{B_1}\right)^p \int_{R_0}^{\infty}\left|e \left(\frac{u}{z}\right) \right|u^{p}\,du.
\end{equation}
In order to conclude, note that the second inequality in~\eqref{eq.Bounds_e_sector}, followed by the first one, and the fact that $e(x)>0$ for $x>0$, together imply that for every $z\in S_{\ga}$ and every $u>0$ we have
$$
|e(u/z)|\le K_3h_{\M}\left(K_4\frac{|z|}{u}\right)\le
\frac{K_3}{K_1}\,e\left(\frac{K_2u}{K_4|z|}\right).
$$
We use this fact, a simple change of variable and the right-hand estimates in~\eqref{eq.boundsMmomentsdc}, and obtain that
\begin{align*}
\int_{0}^{\infty}\left|e \left(\frac{u}{z}\right) \right|u^{p}\,du
&\le
\int_{0}^{\infty}\frac{K_3}{K_1}\,e\left(\frac{K_2u}{K_4|z|}\right) u^{p}\,du\\
&=\frac{K_3}{K_1}\left(\frac{K_4|z|}{K_2}\right)^{p+1}\mo(p)
\le \frac{\mo(0)K_3K_4}{K_1K_2}\left(\frac{K_4B_2}{K_2}\right)^{p}M_p|z|^{p+1}.
\end{align*}
This estimate can be taken into both~\eqref{eq.Estimates_f1} and~\eqref{eq.Estimates_f2}, and from~\eqref{eq.RemainderTwoSums} we easily get that for every $p\in\N_0$,
$$
\left|f(z)-\sum_{n=0}^{p-1}a_nz^n\right|\le \frac{3K_3K_4}{K_1K_2}|\widehat{f}|_{\M,A} \left(\frac{2K_4B_2A}{K_2B_1}\right)^pM_p|z|^p,\quad z\in S_{\ga},
$$
and so $f$ admits $\widehat{f}$ as its uniform $\{\M\}$-asymptotic expansion in $S_{\ga}$. Moreover, recalling the definition~\eqref{eq.NormUniformAsymptFixedType} of the norm in these spaces with uniform asymptotics and fixed type, if we put $c:=2K_4B_2/(K_2B_1)>0$, we see that $f\in\widetilde{\mathcal{A}}^u_{\M,cA}(S_{\ga})$ and
$$
\|f\|_{\M,cA,\overset{\sim}{u}}\le \frac{3K_3K_4}{K_1K_2}|\widehat{f}|_{\M,A},
$$
what proves the continuity of the linear map $T_{\M,A}$.

$(iii)\Rightarrow (iv)$ Immediate for any weight sequence $\M$.

$(iv)\Rightarrow (v)$ It follows from~\eqref{equaContentionSurjectIntervals}, again valid for any weight sequence.

$(v)\Rightarrow (vi)$ This statement is a consequence of Theorem~\ref{teorSurject.dc}.
\end{proof1}

\begin{rema}
The facts in Theorem~\ref{tpral}.$(iii)$ and Proposition~\ref{propcotaderidesaasin}.$(ii)$ together guarantee that for every $\delta\in(0,\ga)$ there exists $c'>0$ such that for every $A>0$ there exists a linear and continuous extension operator from $\C[[z]]_{\M,A}$ into $\mathcal{A}_{\hM,c'A}(S_{\delta})$. In fact, V. Thilliez stated his main result in this regard~\cite[Th.~3.2.1]{Thilliez03} in terms of the existence of such extension operators for every $\delta<\ga(\M)$ and $\M$ a strongly regular sequence.
\end{rema}

The following three corollaries become now clear.

\begin{coro}
Let $\hM$ be a regular sequence,
and $\ga>0$. The following are equivalent:
\begin{itemize}
\item[$(i)$] $\gamma(\M)>\gamma$,

\item[$(ii)$] There exists $\gamma_1>\gamma$ such that the space $\widetilde{\mathcal{A}}^u_{\{\M\}}(S_{\gamma_1})$ contains optimal $\{\M\}$-flat functions.

\item[$(iii)$] There exists $\gamma_1>\gamma$ such  that the Borel map $\widetilde{\mathcal{B}}:\widetilde{\mathcal{A}}^u_{\{\M\}}(S_{\gamma_1})\to \C[[z]]_{\{\M\}}$ is surjective., i.e., $\gamma_1\in\widetilde{S}^u_{\{\M\}}$.
\end{itemize}
\end{coro}

\begin{proof1}
$(ii)\Rightarrow(iii)$ and $(iii)\Rightarrow(i)$ are respectively contained in Theorem~\ref{tpral} and Theorem~\ref{teorSurject.dc}, under weaker hypotheses.
$(i)\Rightarrow(ii)$ is immediately deduced from Proposition~\ref{prop.ExistOptFlatFunct}.
\end{proof1}

As a consequence of \eqref{eq.snqiffgammaMpositive}
and Theorem~\ref{tpral} we get the following result.

\begin{coro}\label{coro.optimalflat}
Let $\hM$ be a regular sequence.
The following are equivalent:
\begin{itemize}
\item[$(i)$] $\M$ satisfies (snq).

\item[$(ii)$] There exists $\gamma>0$ such that the space $\widetilde{\mathcal{A}}^u_{\{\M\}}(S_{\gamma})$ contains optimal $\{\M\}$-flat functions.

\item[$(iii)$] There exists $\gamma>0$ such that the Borel map $\widetilde{\mathcal{B}}:\widetilde{\mathcal{A}}^u_{\{\M\}}(S_{\gamma})\to \C[[z]]_{\{\M\}}$ is surjective. In other words, $\widetilde{S}^u_{\{\M\}}\neq\emptyset$.
\end{itemize}
\end{coro}

Note that, according to Proposition~\ref{propcotaderidesaasin}, in the previous items $(ii)$ and $(iii)$ one could change $\widetilde{\mathcal{A}}^u_{\{\M\}}(S_{\gamma})$ and $\widetilde{S}^u_{\{\M\}}$ into $\mathcal{A}_{\{\hM\}}(S_{\gamma})$ and $S_{\{\hM\}}$, respectively.

\begin{coro}\label{coro.optimalflat0}
Let $\hM$ be a regular sequence,
and $\ga>0$. The following are equivalent:
\begin{itemize}
\item[$(i)$] $\gamma(\M)>\gamma$,

\item[$(ii)$] There exists $\gamma_1>\gamma$ such that the space $\mathcal{A}_{\{\hM\}}(S_{\gamma_1})$ contains optimal $\{\M\}$-flat functions,

\item[$(iii)$] There exists $\gamma_1>\gamma$ such that $\widetilde{\mathcal{B}}:\mathcal{A}_{\{\hM\}}(S_{\gamma_1})\to \C[[z]]_{\{\M\}}$ is surjective, i.e., $\gamma_1\in S_{\{\hM\}}$.
\end{itemize}
\end{coro}

\subsection{Optimal flat functions and strongly regular sequences}

Under the moderate growth condition,
the implication $(ii)\Rightarrow(i)$ in the version of Corollary~\ref{coro.optimalflat} for the space $\mathcal{A}_{\{\hM\}}(S_{\gamma})$
can be shown independently by using a result from J. Bruna~\cite{Brunaext80}, where a precise formula for nontrivial flat functions in Carleman-Roumieu ultradifferentiable classes, appearing in a work of T. Bang~\cite{Bang53}, is exploited. For the sake of completeness, we will present this proof below.

\begin{theo}\label{teor.optimalflat1}
Let $\M$ be a weight sequence satisfying (mg). If there exists $\gamma>0$ such that $\mathcal{A}_{\{\hM\}}(S_{\gamma})$ contains optimal $\{\M\}$-flat functions, then $\M$ is strongly regular.
\end{theo}

The proof requires two auxiliary results which we state and prove now.

First, given a weight sequence $\M$, the sequence of quotients $\m=(m_p)_{p\in\N_0}$ is nondecreasing and tends to infinity, but it can happen that it remains constant on large intervals $[p_0,p_1]$ of indices, so that the counting function $\nu_{\m}$ defined in \eqref{countinfunctionnu} yields $\nu_{\m}(m_{p_0})=\nu_{\m}(m_{p_1})=p_1+1$. However, in some applications or proofs it would be convenient to have $\nu_{\m}(m_p)=p+1$ for all $p\ge 0$. This can be assumed without loss of generality by the following result.

\begin{lemma}\label{strictincreasinglemma}
Let $a=(a_p)_{p\ge 1}$ be a nondecreasing sequence of positive real numbers satisfying $\lim_{p\rightarrow+\infty}a_p=+\infty$ (it suffices that $a_{p-1}<a_{p}$ holds true for infinitely many indices $p$). Then there exists a sequence $b=(b_p)_{p\ge 1}$ of positive real numbers such that $p\mapsto b_p$ is strictly increasing and satisfies
$$
0<\inf_{p\ge 1}\frac{b_p}{a_p}\le\sup_{p\ge 1}\frac{b_p}{a_p}<+\infty.
$$
\end{lemma}

So, in the language of weight sequences, we prove that for any weight sequence $\M$ there exists a strongly equivalent weight sequence $\L$ (and so $\M\approx\L$) such that $\nu_{\boldsymbol{\ell}}(\ell_p)=p+1$ for all $p\in\N_0$. Note that equivalent weight sequences define the same Carleman-Roumieu ultraholomorphic classes and associated weighted classes of formal power series.

\begin{proof1}
Since $a$ is nondecreasing and $\lim_{p\rightarrow+\infty}a_p=+\infty$ there exists a sequence $(p_j)_{j\ge 1}$ of indices such that  $a_{p_j-1}<a_{p_j}=\dots=a_{p_{j+1}-1}<a_{p_{j+1}}$ for all $j\ge 1$ (and so $p_1\ge 2$). For all $j\ge 1$ we have now $a_{p_j}/(a_{p_j-1})>1+\varepsilon_j$ for a sequence $(\varepsilon_j)_{j\ge 1}$ with possibly small strictly positive numbers $\varepsilon_j$. Finally we put $p_0:=1$.\vspace{6pt}

We take some arbitrary $A>1$ and choose $\delta_j>0$ small enough so as to have $(1+\delta_j)^{p_{j+1}-p_j-1}\le\min\{A,1+\varepsilon_{j+1}\}$. Then the sequence $(\delta_j)_{j\ge 0}$ satisfies
\begin{equation}\label{strict1}
(1+\delta_j)^{p_{j+1}-p_j-1}\le 1+\varepsilon_{j+1}<\frac{a_{p_{j+1}}}{a_{p_{j+1}-1}},\quad (1+\delta_j)^{p_{j+1}-p_j-1}\le A,\quad j\ge 0.
\end{equation}
We define now $b$ as follows:
\begin{equation}\label{strict2}
b_q:=a_q\;\text{if}\;q=p_j,\;j\ge 0,\hspace{20pt}b_q:=(1+\delta_j)b_{q-1}\;\text{if}\;1+p_j\le q\le p_{j+1}-1,\;j\ge 0.
\end{equation}

So we have by iteration $b_q=(1+\delta_j)^{q-p_j}b_{p_j}=(1+\delta_j)^{q-p_j}a_{p_j}=(1+\delta_j)^{q-p_j}a_q>a_q$  for all $q$ with $1+p_j\le q\le p_{j+1}-1$, $j\ge 0$. On each such interval of indices the mapping $q\mapsto b_q$ is now clearly strictly increasing since $1+\delta_j>1$ for all $j$.
Moreover, by the first half in \eqref{strict1}, we have $b_{p_{j+1}-1}=(1+\delta_j)^{p_{j+1}-p_j-1}a_{p_j}<b_{p_{j+1}}$.
Hence the sequence $q\mapsto b_q$ is strictly increasing.\vspace{6pt}

By definition \eqref{strict2} we have $b_q=a_q$ for all $q=p_j$, $j\ge 0$, and $b_q>a_q$ otherwise. We conclude if we show that $b_q\le A a_q$ for all $q$ with $1+p_j\le q\le p_{j+1}-1$, $j\ge 0$. For this, since $q\mapsto b_q$ is strictly increasing, it suffices to observe that, thanks to the second half in \eqref{strict1}, we have $b_{p_{j+1}-1}=(1+\delta_j)^{p_{j+1}-p_j-1}a_{p_j}\le Aa_{p_j}=Aa_{p_{j+1}-1}$.
\end{proof1}

The second result is the following.

\begin{lemma}\label{equivalenceofdual}
Let $\M$ be a weight sequence. Then $\M$ satisfies (mg) if and only if
$\omega_\M(t)=O(\nu_{\m}(t))$ as $t\rightarrow+\infty$.
\end{lemma}

\begin{proof1}
The condition (mg) for $\M$ is equivalent to $m_n\le A(M_n)^{1/n}$ for some $A\ge 1$ and all $n\in\N$ (e.g., see~\cite[Lemma 2.2]{RainerSchindlExtension17}).
It is also known that $\omega_\M(m_n)=\log\left(m_{n}^n/M_n\right)$ for $n\in\N$ (see~\cite[Chapitre I]{Mandelbrojt}). So, if $m_{n-1}\le t<m_{n}$ for some $n\ge 1$, we get
$$
\omega_\M(t)\le\omega_\M(m_n)=n\log\left(\frac{m_n}{M_n^{1/n}}\right)\le n\log(A)=\log(A)\nu_{\m}(t),
$$
that is, $\omega_\M(t)=O(\nu_{\m}(t))$ as $t\rightarrow+\infty$.\vspace{6pt}

Conversely, suppose that there exists $A\ge 1$ such that $\omega_\M(t)\le A\nu_{\m}(t)$ for all $t\ge m_0$. By \cite[Lemma 2.2]{RainerSchindlExtension17}, (mg) for $\M$ holds true if and only if there exists $H\geq 1$ such that for all $t$ large enough one has $2\nu_{\m}(t)\leq\nu_{\m}(H t)+H$, and this we will prove. Take $H\ge \exp(2A)$ and $t\ge m_0$. Using \eqref{assofuncintegral}, and since $\nu_{\m}$ is nondecreasing, we estimate
\begin{align*}
\nu_{\m}(Ht)&\ge A^{-1}\omega_\M(Ht)= A^{-1}\int_{m_0}^{Ht}\frac{\nu_{\m}(\lambda)}{\lambda}d\lambda\ge A^{-1}\int_t^{Ht}\frac{\nu_{\m}(\lambda)}{\lambda}d\lambda
\\&
\ge A^{-1}\nu_{\m}(t)\int_t^{Ht}\frac{1}{\lambda}d\lambda= A^{-1}\log(H)\nu_{\m}(t)\ge 2\nu_{\m}(t),
\end{align*}
as desired.

We mention that an alternative, more abstract proof can be based in the theory of O-regular variation and Matuszewska indices for  functions. By \cite[Th. 4.4]{JimenezSanzSchindlIndex} we have that the lower Matuszewska indices of $\nu_{\m}$ and $\omega_{\M}$ agree, that is, $\beta(\nu_{\m})= \beta(\omega_\M)$, and by \cite[Cor.2.17 and Cor. 4.2]{JimenezSanzSchindlIndex} we know $\M$ has (mg) if and only if $\beta(\nu_{\m})>0$. So, if $\beta(\nu_{\m})>0$, by \cite[Th. 4.3]{JimenezSanzSchindlIndex} we have that $\liminf_{t\to\infty} \frac{\nu_{\m}(t)}{\omega_\M(t)}>0$, and we deduce that $\omega_\M(t)=O(\nu_{\m}(t))$ as $t\rightarrow+\infty$. Conversely, if $\omega_\M(t)=O(\nu_{\m}(t))$ as $t\rightarrow+\infty$, then  $\liminf_{t\to\infty} \frac{\nu_{\m}(t)}{\omega_\M(t)}>0$, so by \cite[Th. 4.3]{JimenezSanzSchindlIndex} we have that   $\beta(\omega_{\M})>0$, and we are done.
\end{proof1}

\noindent\textbf{Proof of Theorem~\ref{teor.optimalflat1}. } We follow the proof of necessity for \cite[Th.~2.2]{Brunaext80}. By Lemma \ref{strictincreasinglemma} and the remark following it, we can assume without loss of generality that $\m$ is strictly increasing.

Let $G$ be an optimal $\{\M\}$-flat function in $\mathcal{A}_{\{\hM\}}(S_{\gamma})$ for some $\gamma>0$.
So, there exists some $A>0$ such that
$$
p_{\M,A}(G):=\sup_{n\in\N_0,x\in(0,+\infty)}\frac{|G^{(n)}(x)|}{A^n n!M_n}<+\infty.
$$
This shows that the Carleman-Roumieu ultradifferentiable class $\mathcal{E}_{\{\hM\}}((-\varepsilon,+\infty))$, consisting of all smooth complex-valued functions $g$ defined on the interval $(-\varepsilon,\infty)$ for some $\varepsilon>0$, and such that $$\sup_{n\in\N_0,x\in(-\varepsilon,+\infty)}\frac{|g^{(n)}(x)|}{D^n n!M_n}<+\infty
$$
for suitable $D>0$, contains nontrivial flat functions (it suffices to extend $G$ by 0 for $x\in(-\varepsilon,0]$). Then, the well-known Denjoy-Carleman theorem (e.g., see \cite[Th.~1.3.8]{hoermander}) yields that $\M$ satisfies (nq).

Let now
$$
R_n:=\sum_{k\ge n}\frac{1}{(k+1)m_k}<+\infty,\quad n\in\N_0,
$$
and let the function $\overline{h}$ be defined by $\overline{h}(t):=n$ if $R_{n+1}<t\le R_n$, $n\in\N_0$.

By \cite[$(14)$, p. 142]{Bang53} we obtain that
$$
G(x)=|G(x)|\le p_{\M,A}(G)\exp\left(-\overline{h}(Aex)\right),\quad x\in(0,+\infty).
$$
Combining this with \eqref{optimalflatleft}, with~\eqref{functionhequ2} and setting $C:=p_{\M,A}(G)$, we get
$$
\exp\left(\overline{h}(Aex)\right)\le\frac{C}{G(x)}\le C K_1^{-1}\exp(\omega_{\M}(1/(K_2x))),\quad x>0.
$$
If we put $t=Aex$ and $B:=Ae/K_2$, we obtain that for every $t>0$,
\begin{equation}\label{optimalflattheoremequ1}
\overline{h}(t)\le \log(C K_1^{-1}) +\omega_{\M}(B/t).
\end{equation}
By Lemma \ref{equivalenceofdual}, there exists $C_1\ge 1$ such that $\omega_\M(s)\le C_1\nu_{\m}(s)+C_1$ for $s>0$. Choosing $t=B/m_n$ in~\eqref{optimalflattheoremequ1}, we see that
$$
\overline{h}(B/m_n)\le \log(C K_1^{-1}) +\omega_{\M}(m_n)\le \log(C K_1^{-1}) +C_1\nu_{\m}(m_n)+C_1=\log(C K_1^{-1}) +C_1(n+1)+C_1,
$$
since $\m$ is strictly increasing. Hence, $\overline{h}(B/m_n)\le C_2(n+1)$ for some $C_2\in\N$ and all $n\in\N_0$. By definition of $\overline{h}$, we get $R_{C_2(n+1)+1}\leq B/m_n$, i.e.,
$$
m_n\sum_{k\ge C_2(n+1)+1}\frac{1}{(k+1)m_k}\le B,\quad n\in\N_0.
$$
Finally,
\begin{align*}
	m_n\sum_{k\ge n}\frac{1}{(k+1)m_k}&=m_n\sum_{k\ge C_2(n+1)+1}\frac{1}{(k+1)m_k}+m_n\sum_{k=n}^{C_2(n+1)}\frac{1}{(k+1)m_k}\\
	&\le B+m_n\frac{n(C_2-1)+C_2+1}{(n+1)m_{n}}\le B+2C_2,
\end{align*}
which is (snq) for $\M$.\hfill$\Box$

\section{Construction of optimal flat functions for a family of non strongly regular sequences}\label{sectConstrOptFlatNon_SR_Seq}

As deduced in Theorem~\ref{tpral}, the construction of optimal $\{\M\}$-flat functions in sectors within an ultraholomorphic class, given by a regular sequence $\hM$, provides extension operators and surjectivity results.
Although such general construction has been shown in Proposition~\ref{prop.ExistOptFlatFunct}, we wish to present here a family of (non strongly) regular sequences for which an alternative, more explicit technique works.

We recall that, for logarithmically convex sequences $(M_p)_{p\in\N_0}$, the condition  (dc) is equivalent to the condition $\log(M_p)=O(p^2)$, $p\to\infty$ (see~\cite[Ch.~6]{Mandelbrojt}). On the other hand, the condition (mg) implies that the sequence is below some Gevrey order (there exists $\a>0$ such that $M_p=O(p!^{\a})$ as $p\to\infty$; see e.g.~\cite{Matsumoto84,Thilliez03}).

We will work, for $q>1$ and $1<\sigma\le 2$, with the sequences $\M_{q,\sigma}:=(q^{p^\sigma})_{p\in\N_0}$. They are clearly weight sequences and, by~\eqref{equa.indice.gammaM.casicrec}, it is immediate that $\ga(\M_{q,\sigma})=\infty$, so they satisfy (snq) (see~\eqref{eq.snqiffgammaMpositive}). According to the previous comments, they satisfy (dc) but not (mg). So, $\hM_{q,\sigma}$ is regular, but $\M_{q,\sigma}$ is not strongly regular.

The case $\sigma=2$ is well-known, as it corresponds
to the so-called \emph{$q$-Gevrey sequences}, appearing in the study of formal solutions for $q$-difference equations.

First, we will construct a holomorphic function on $\C\setminus(-\infty,0]$ which will provide, by restriction, an optimal $\{\M_{q,\sigma}\}$-flat function in any unbounded sector $S_{\ga}$ with $0<\ga<2$. Subsequently, we will obtain such functions on general sectors of the Riemann surface $\mathcal{R}$ of the logarithm by ramification. This, according to Theorems~\ref{teorSurject.dc} and~\ref{tpral}, agrees with the fact that $\ga(\M_{q,\sigma})=\infty$.

\subsection{Flatness in the class given by $\Ms$}

It will be convenient to note that for a fixed $\sigma\in(1,2]$, there exists a unique $s\ge 2$ such that $\sigma=s/(s-1)$.

We start by suitably estimating the function
$$\o_{\Ms}(t)= \sup_{p\in\mathbb{N}_0} \ln\left(\frac{t^{p}}{q^{p^\sigma}}\right)=\sup_{p\in\mathbb{N}_0} (p\ln(t)-p^{s/(s-1)}\ln(q)),\quad t>0.
$$
Due to the fact that $\o_{\Ms}(t)=0$ for $t\leq1$ (since $m_0=M_1/M_0=M_1=q>1$ and by~\eqref{assofuncintegral}), we will restrict our attention to the case $t>1$.
Obviously, $\o_{\Ms}(t)$ is bounded above by the supremum of $x\ln(t)-x^{s/(s-1)}\ln(q)$ when $x$ runs over $(0,\infty)$, which is easily obtained by elementary calculus and occurs at the point
$$
x_0=\left(\frac{(s-1)\ln(t)}{s\ln(q)}\right)^{s-1}.
$$
If we put
\begin{equation}\label{eq.def.b_qs}
b_{q,s}:=\frac{1}{s}\left(\frac{s-1}{s\ln(q)}\right)^{s-1},
\end{equation}
then
\begin{equation}\label{o_Ms_upbound}
	\o_{\Ms}(t)\leq \left(\frac{(s-1)\ln(t)}{s\ln(q)}\right)^{s-1}\ln(t)- \left(\frac{(s-1)\ln(t)}{s\ln(q)}\right)^{s}\ln(q)=b_{q,s} \ln^s(t),\quad t>1.
\end{equation}
On the other hand, for $t> q^{s/(s-1)} $ (what amounts to $x_0>1$) we also have that $\o_{\Ms}(t)$ is at least the value of $x\ln(t)-x^{s/(s-1)}\ln(q)$ at $x=\flo{x_0}$, that is,
\begin{align}
\o_{\Ms}(t)&\ge \flo{\left(\frac{(s-1)\ln(t)}{s\ln(q)}\right)^{s-1}}\ln(t)- \flo{\left(\frac{(s-1)\ln(t)}{s\ln(q)}\right)^{s-1}}^{s/(s-1)} \ln(q)\nonumber\\
&\ge \left(\left(\frac{(s-1)\ln(t)}{s\ln(q)}\right)^{s-1}-1\right)\ln(t)- \left(\frac{(s-1)\ln(t)}{s\ln(q)}\right)^s\ln(q)\nonumber\\
&=b_{q,s} \ln^s(t)-\ln(t).\label{o_Ms_interbound}
\end{align}

\begin{lemma}
For every $t\ge q^{2s/(s-1)}$ it holds
\begin{equation}\label{eq.EstimateOmegaQSigma}
b_{q,s}\ln^s(t)-\ln(t)\ge b_{q,s}\ln^s\left(\frac{t}{q^{s/(s-1)}}\right)- \ln\left(q^{s/(s-1)}\right).
\end{equation}
\end{lemma}
\begin{proof1}
Observe that every $t\ge q^{2s/(s-1)}$ may be written as $t=q^{ys/(s-1)}$ for some $y\ge 2$. Then, we have that
\begin{equation*}
b_{q,s}\ln^s(t)-b_{q,s}\ln^s\left(\frac{t}{q^{s/(s-1)}}\right)=
b_{q,s}\left(\frac{s\ln(q)}{s-1}\right)^s\big(y^s-(y-1)^s\big)=
\frac{\ln(q)}{s-1}\big(y^s-(y-1)^s\big).
\end{equation*}
By the mean value theorem, $y^s-(y-1)^s>s(y-1)^{s-1}$, and since $s\ge 2$ and $y\ge 2$, we have $(y-1)^{s-1}\ge y-1$. So we deduce that
$$
\frac{\ln(q)}{s-1}\big(y^s-(y-1)^s\big)>\frac{s\ln(q)}{s-1}(y-1)=
\ln(t)-\ln\left(q^{s/(s-1)}\right),
$$
as desired.
\end{proof1}

Combining~\eqref{o_Ms_upbound} with~\eqref{o_Ms_interbound} and~\eqref{eq.EstimateOmegaQSigma}, and using~\eqref{functionhequ2}, we get
\begin{equation}\label{Ms_bound_h}
\exp\left(-b_{q,s}\ln^s\left(\frac{1}{t}\right)\right)	\leq h_{\Ms}(t) \leq q^{s/(s-1)}\exp\left(-b_{q,s}\ln^s\left(\frac{1}{q^{s/(s-1)}t}\right)\right), \qquad 0<t\leq q^{-2s/(s-1)}.
\end{equation}

We can say that these estimates express optimal $\{\Ms\}$-flatness.

\subsection{Optimal $\{\Ms\}$-flat function in $S_2$}

The estimates in~\eqref{Ms_bound_h} suggest considering the function $\exp\big(-b_{q,s}\log^s\big(1/z\big)\big)$, with, say, principal branches, as a candidate for providing optimal flat functions. However, its analyticity in wide sectors is not guaranteed. Moreover, even in small sectors around the direction $d=0$, its behaviour at $\infty$ might not be as desired: For example, when $s=2$ it tends to 0 as $0<x\to\infty$, what excludes the possibility of proving the
inequality in~\eqref{optimalflatleft}.

Because of these reasons, we will first define a suitably modified function in the sector $S_2=\C\setminus(-\infty,0]$, prove its flatness there, and then turn to general sectors by composing it with an appropriate ramification.

We define
\begin{equation}\label{eq.defG2}
G_2^{q,s}(z):=\exp\left(-b_{q,s}\log^s\left(1+\frac{1}{z}\right)\right),\quad z\in S_2,
\end{equation}
where the principal branch of the logarithm is chosen for both $\log$ and the power $w\mapsto w^s=\exp(s\log(w))$ involved.
Observe that if $z\in S_2$, then $1+1/z \in\mathbb{C}\setminus(-\infty,1]$, and so
$\log(1+1/z)=\ln(|1+1/z|)+i\arg(1+1/z)\notin(-\infty,0]$. This ensures that the map
\begin{align*}
z&\mapsto \log^{s}\left(1+\frac{1}{z}\right)= \exp\left(s\log\left(\log(1+\frac{1}{z})\right)\right)
\end{align*}
is also holomorphic in $S_2$, and so is $G_2^{q,s}$.

In order that $G_2^{q,s}$ is an optimal $\{\Ms\}$-flat function in $S_2$, we are only left with proving the
estimates~\eqref{optimalflatleft} and \eqref{optimalflatright}.
It turns out to be more convenient to work with the associated kernel
$$
e_2(z):=G_2^{q,s}(1/z)=\exp(-b_{q,s}\log^s(1+z)),\quad z\in S_2,
$$
and verify the following result.

\begin{lemma}\label{e2-bounds}
There exist positive constants $C_1,C_2$ such that
$$
|e_2(z)|\leq C_1 e_2(C_2 |z|),\quad z\in S_2.
$$
\end{lemma}
\begin{proof1}
%
In the first place, we observe that for every $z\in S_2$,
\begin{equation}\label{eq.RealPartPowerLog}
\Re(\log^s(z+1))= |\log^s(z+1)|\cos(\arg(\log^s(z+1)))=
|\log(z+1)|^s\cos(s\arg(\log(z+1))).
\end{equation}
Now,
\begin{equation}\label{arglog-inequality} s|\arg(\log(z+1))|= s\left|\arctan\left(\frac{\arg(z+1)}{\ln|z+1|}\right)\right|\leq s\left|\arctan\left(\frac{\pi}{\ln|z+1|}\right)\right|.
\end{equation}
Hence, setting
$$
R_0:=1+\exp\left(\frac{\pi}{\tan\left(\pi/(2s)\right)}\right)\ge 2,
$$
we get that $|z|>R_0$ implies that $|z+1|>R_0-1\ge 1$, and therefore $\ln|z+1|> 0$ and
$$
\frac{\pi}{\ln|z+1|}<\tan\left(\frac{\pi}{2s}\right).
$$
From this and \eqref{arglog-inequality} we deduce that $\cos(s\arg(\log(z+1)))>0$. Then, continuing with~\eqref{eq.RealPartPowerLog},
\begin{align}\label{realpart_logs_R0}
\Re(\log^s(z+1))&\ge |\Re(\log(z+1))|^s\cos(s\arg(\log(z+1)))
\nonumber\\
&=\ln^s|z+1|- \ln^s|z+1|\frac{\sin^2(s\arg(\log(z+1)))}{1+\cos(s\arg(\log(z+1)))}.	 
\end{align}
Now, from the equality in~\eqref{arglog-inequality} we see that $s\arg(\log(z+1))\to 0$ as $z\to \infty$ in $S_2$, and moreover
$$
\lim_{\substack{z\to \infty\\z\in S_2}}\left[ \left(\frac{\sin^2(s\arg(\log(z+1)))}{1+\cos(s\arg(\log(z+1)))} \right) \Big/\left(\frac{s^2\arg^2(z+1)}{2\ln^2|z+1|}\right)\right]=1.
$$
Therefore, there exist $R_1\geq R_0$ and $C>0$ such that
	\begin{equation*}
		 \frac{\sin^2(s\arg(\log(z+1)))}{1+\cos(s\arg(\log(z+1)))}\leq C \frac{1}{\ln^2|z+1|},\qquad |z|>R_1.
	\end{equation*}
We deduce from \eqref{realpart_logs_R0} that for $z\in S_2$ with $|z|>R_1$,
	\begin{equation}\label{realpart_logs_R1}
		\Re(\log^s(z+1))\geq \ln^s|z+1|-C\ln^{s-2}|z+1| \geq  \ln^s(|z|-1)-C\ln^{s-2}(|z|+1).
	\end{equation}
We would be almost done if we obtain, for the right-hand side in~\eqref{realpart_logs_R1}, a lower bound in terms of, say, $\ln^s(1+|z|/2)$ for $|z|$ sufficiently large.

This is easy in case $s=2$, for it suffices to take $|z|>4$ in order to have $3<1+|z|/2<|z|-1$, and so if $|z|\ge R_2:=\max\{R_1,4\}$ we have
$$
\Re(\log^s(z+1))\geq \ln^s(|z|-1)-C\ge \ln^s\left(1+\frac{|z|}{2}\right)-C.
$$
In case $s>2$, it is not difficult to check that
$$
\lim_{x\to+\infty}\left(\ln^s(x-1)-C\ln^{s-2}(x+1)- \ln^s\left(1+\frac{x}{2}\right)\right)=+\infty,
$$
so that, according to~\eqref{realpart_logs_R1}, there exists $R_2\ge R_1$ such that for $z\in S_2$ with $|z|\ge R_2$ one has
\begin{equation*}
\Re(\log^s(z+1))\geq \ln^s\left(1+\frac{|z|}{2}\right).
\end{equation*}
In any case, we can deduce an upper estimate of the form
\begin{align*}
|e_2(z)|&=\exp\big(-b_{q,s}\Re(\log^s(z+1))\big)\\
&\leq  e^C\exp\left(-b_{q,s}\ln^s\left(1+\frac{|z|}{2}\right)\right)= e^Ce_2\left(\frac{|z|}{2}\right), \quad z\in S_2,\,|z|> R_2.
\end{align*}
Finally, since the function $|e_2(z)|$ stays bounded and bounded away from 0 for bounded $|z|$ (in particular, it tends to 1 when $z$ tends to 0 in $S_2$), the previous estimate can be extended to the whole of $S_2$ by suitably enlarging the constant $C$.
\end{proof1}

We are ready for the main objective of this section.

\begin{theo}\label{teor.FlatFunctionS2}
The function $G_2^{q,s}$ defined in~\eqref{eq.defG2} is an optimal $\{\Ms\}$-flat function in $S_2$.
\end{theo}

\begin{proof1}
The previous lemma ensures that there exist positive constants $C_1$, $C_2$ such that
\begin{equation}\label{G2-bound}
|G_2^{q,s}(z)|\leq C_1 \exp\left(-b_{q,s}\ln^s\left(1+\frac{C_2}{|z|}\right)\right),\quad z\in S_2.
\end{equation}
Observe that this inequality guarantees that $|G_2^{q,s}(z)|$ is bounded.
As the same is true for $h_{\Ms}(t)$ for every $t\ge t_0$ and any fixed $t_0>0$ (see Lemma~\ref{functionhproperties}), we only need to check the estimate~\eqref{optimalflatright} for small enough $|z|$.

For $|z|\le C_2$ it is clear that $\ln(1+C_2/|z|)>\ln(C_2/|z|)\ge 0$. Then, from
~\eqref{Ms_bound_h} we have that
\begin{align*}
|G_2^{q,s}(z)|&\leq C_1\exp\left(-b_{q,s}\ln^s\left(1+\frac{C_2}{|z|}\right)\right)\\
&\leq C_1 \exp\left(-b_{q,s}\ln^s\left(\frac{C_2}{|z|}\right)\right)
\leq C_1 h_{\Ms}\left(\frac{|z|}{C_2}\right),\quad |z|\leq C_2q^{-2s/(s-1)},
\end{align*}
and we have proved
~\eqref{optimalflatright}.

Now, let us note that $G_2^{q,s}(x)$ is bounded away from 0 as soon as $x\ge r$ for any fixed $r>0$, since then
	$$
	\exp\left(-b_{q,s}\ln^s\left(1+1/r\right)\right)\leq G_2^{q,s}(x).
	$$
Again, we only need to check the estimate~\eqref{optimalflatleft} for small enough $x$. Indeed, we have for $x>0$ that
\begin{equation*}
		G_2^{q,s}(x)= \exp\left(-b_{q,s}\ln^s\left(\frac{1}{x}\right)\right) \exp\left(-b_{q,s}\left[\ln^s\left(1+\frac{1}{x}\right)- \ln^s\left(\frac{1}{x}\right)\right]\right).
\end{equation*}
The mean value theorem gives that $\ln^s(1+1/x)-\ln^s(1/x)$ tends to zero if $x\searrow 0$, and we deduce that there exists $L$ such that
\begin{equation*}
		G_2^{q,s}(x)\geq  L\exp\left(-b_{q,s}\ln^s\left(\frac{1}{x}\right)\right),
		\quad x\leq {\color{blue}q^{-s/(s-1)}}.
\end{equation*}
The second inequality in~\eqref{Ms_bound_h} implies now that, as long as $x\le q^{-s/(s-1)}$, we have
\begin{equation*}
		G_2^{q,s}(x)\geq  Lq^{-s/(s-1)} h_{\Ms}\left(\frac{x}{q^{s/(s-1)}}\right),
\end{equation*}
and so \eqref{optimalflatleft} holds.
%
%
\end{proof1}

\subsection{Optimal $\{\Ms\}$-flat function in arbitrary sectors}

Let us consider a sector $S_\gamma\subset\mathcal{R}$ with $\gamma>2$, and define the function
\begin{equation}\label{eq.defGgamma}
G_{\ga}^{q,s}(z):=\exp\left(-b_{q,s}\left(\frac{\gamma}{2}\right)^s \log^s\left(1+z^{-2/\gamma}\right)\right) =\left(G_2^{q,s}(z^{2/\gamma})\right)^{(\gamma/2)^s},\quad z\in S_{\ga}.
\end{equation}
The map $z\mapsto z^{2/\gamma}$ is holomorphic from $S_{\ga}$ into $S_2$, and so $G_{\ga}^{q,s}$ is holomorphic in $S_{\ga}$. We will prove that this function is an optimal $\{\Ms\}$-flat function in this sector.

As before, we consider the kernel
$$e_\gamma(z):=G_{\ga}^{q,s}(1/z)= \exp\left(-b_{q,s}\left(\frac{\gamma}{2}\right)^s \log^s\left(1+z^{2/\gamma}\right)\right) =\left(e_2(z^{2/\gamma})\right)^{(\gamma/2)^s},\quad z\in S_{\ga}.$$

\begin{lemma}\label{egamma-bounds}
There exist constants $B_1,B_2>0$ such that
\begin{equation}\label{eq.Bounds_e_gamma}
|e_\gamma(z)|\leq B_1 e_2(B_2 |z|),\quad z\in S_\gamma.
\end{equation}
\end{lemma}
\begin{proof1}
According to the definition of $e_\ga$ and by applying Lemma~\ref{e2-bounds}, there exist constants
$C_1,C_2>0$ such that for every $z\in S_\ga$ one has
$$
|e_\gamma(z)|=\left|e_2(z^{2/\gamma})\right|^{(\gamma/2)^s}\leq \left(C_1 e_2(C_2 |z|^{2/\gamma})\right)^{(\ga/2)^s}.$$
We recall that the function $|e_2(z)|$ stays bounded
for bounded $|z|$; from the previous estimates, the same can be said about $|e_{\ga}(z)|$, and so we can prove~\eqref{eq.Bounds_e_gamma} by restricting our considerations to large enough values of $|z|$ and well chosen $B_2>0$, and then suitably enlarging the constant $B_1>0$ involved.
Let us observe that
\begin{align*}
\left(e_2(C_2 |z|^{2/\gamma})\right)^{(\ga/2)^s}&= \exp\left(-b_{q,s}\ln^s\left[(1+C_2|z|^{2/\ga})^{\ga/2}\right]\right),\\ e_2(B_2|z|)&=\exp\left(-b_{q,s}\ln^s(1+B_2|z|)\right).
\end{align*}
So, we will be done
if we see that
$$
\ln^s(1+B_2|z|)-\ln^s[(1+C_2|z|^{2/\ga})^{\ga/2}]
$$
admits an upper bound for large enough $|z|$ and suitably chosen $B_2>0$.
But this follows from the clear fact that
\begin{equation*}
\ln^s(1+B_2|z|)-\ln^s\left[\left(1+C_2|z|^{2/\ga}\right)^{\ga/2}\right]
\sim -s\ln\left(\frac{C_2^{\ga/2}}{B_2}\right)\ln^{s-1}(1+B_2|z|), \quad |z|\to\infty,
\end{equation*}
where $\sim$ means that the quotient of both expressions tends to 1.
Indeed, in view of this equivalence it suffices to choose any $B_2<C_2^{\ga/2}$ in order to have the desired estimation for suitably large $B_1$ and $|z|$.

\end{proof1}

\begin{coro}\label{coro.FlatFunctionSgamma}
The function $G_{\ga}^{q,s}$ defined in~\eqref{eq.defGgamma} is an optimal $\{\Ms\}$-flat function in $S_{\ga}$.
\end{coro}

\begin{proof1}
By the previous lemma, there exist $B_1,B_2>0$ such that
\begin{equation*}
|G_{\ga}^{q,s}(z)|\leq B_1 \exp\left(-b_{q,s}\ln^s\left(1+\frac{B_2}{|z|}\right)\right),\quad  z\in S_\ga.
\end{equation*}
Note that this estimate is essentially that in~\eqref{G2-bound}, and so the conclusion follows in exactly the same way as in the proof of Theorem~\ref{teor.FlatFunctionS2}.
\end{proof1}

\begin{rema}
We mention that a similar approach has been followed in the preprint~\cite{JimenezLastraSanzRapidGrowth}, by A. Lastra and the first and third authors, in order to construct extension operators for the ultraholomorphic classes associated with the sequences $\M^{\tau,\sigma}=(p^{\tau p^{\sigma}})_{p\in\N_0}$, for $\tau>0$ and $\sigma\in(1,2)$. These sequences have appeared in a series of papers by S. Pilipovi{\'c}, N. Teofanov and F. Tomi{\'c}~\cite{ptt15,ptt16,ptt,ptt21}, inducing ultradifferentiable spaces of so-called extended Gevrey regularity. However, in that case the construction of suitable kernels for our technique involves the Lambert function, whose handling is not so convenient. This fact has caused our results to be available only in sectors strictly contained in $S_2$, in spite of the fact that $\ga(\M_{\tau,\sigma})=\infty$, what would in principle allow for such extension operators to exist in sectors of arbitrary opening.
\end{rema}

\section{Convolved sequences, flat functions and extension results}

We show in this section that whenever two weight sequences are given and there exist optimal flat functions in the respectively associated classes, then optimal flat functions exist in the class defined by the so-called convolved sequence as well (given by the point-wise product). Moreover, the extension technique works if one of the convolved sequences satisfies (dc).\vspace{6pt}

On the one hand the abstract statement is a straight-forward consequence of a result by H. Komatsu, see Remark \ref{convolvesequrem0} for more details. On the other hand this approach can be useful for constructing (counter-)examples. In general even for nice sequences the convolved sequence can behave complicated, see Sect. \ref{convolveexamplesect}, and so a direct explicit construction of optimal flat functions in the class defined by the convolved sequence will be challenging.

\subsection{Convolved sequences}

Let $\mathbb{M}^1=(M_p^1)_{p\in\N_0}$, $\mathbb{M}^2=(M_p^2)_{p\in\N_0}$ be two sequences of positive real numbers, then the \emph{convolved sequence} $\L:=\mathbb{M}^1\star\mathbb{M}^2$ is $(L_p)_{p\in\N_0}$ given by
\begin{equation*}
L_p:=\min_{0\le q\le p}M^1_qM^2_{p-q},\quad p\in\mathbb{N}_0,
\end{equation*}
see \cite[$(3.15)$]{komatsu}. Hence, obviously $\mathbb{M}^1\star\mathbb{M}^2=\mathbb{M}^2\star\mathbb{M}^1$.

For all $p\in\mathbb{N}_0$ we have $L_p\le\min\{M^1_0M^2_p,M^2_0M^1_p\}$. So, if in addition $M^1_0=M^2_0=1$, then we get $L_0=1$ and
\begin{equation}\label{convolvedbelow}
    L_p\le\min\{M^1_p,M^2_p\}, \quad p\in\N_0.
\end{equation}

Given $\M=(M_p)_{p\in\N_0}$ with $M_0=1$, put $\L=(L_p)_{p\in\N_0}=\M\star\M$. The condition (mg) states precisely that there exists $A>0$ such that $M_p\le A^p L_p$ for every $p\in\N_0$; according to~\eqref{convolvedbelow}, $\M$ satisfies (mg) if and only if $\M$ and $\M\star\M$ are equivalent.

\begin{rema}\label{convolvesequrem0}
Let $\M,\mathbb{M}^1,\mathbb{M}^2$ be weight sequences.
\begin{itemize}
\item[(i)] In \cite[Lemma 3.5]{komatsu} the following facts are shown: $\mathbb{M}^1\star\mathbb{M}^2$ is again a weight sequence. The corresponding quotient sequence $\bm^1\star\bm^2$ is obtained when rearranging resp. ordering the sequences $\bm^1$ and $\bm^2$ in the order of growth. This yields, by definition of the counting function (see~(6)), that for all $t\ge 0$ one has
$$\nu_{\bm^1\star\bm^2}(t)= \nu_{\bm^1}(t)+\nu_{\bm^2}(t);$$
so, by~(7) we get
$$\omega_{\mathbb{M}^1\star\mathbb{M}^2}(t)= \omega_{\mathbb{M}^1}(t)+\omega_{\mathbb{M}^2}(t),\quad t\ge 0,$$
and by~(4) we obtain
\begin{equation}\label{functionhforconvolved}
h_{\mathbb{M}^1\star\mathbb{M}^2}(t)= h_{\mathbb{M}^1}(t)h_{\mathbb{M}^2}(t),\quad t>0.
\end{equation}
\item[(ii)] If either $\mathbb{M}^1$ or $\mathbb{M}^2$ has (dc), then $\mathbb{M}^1\star\mathbb{M}^2$ as well: As said before, for sequences $(M_p)_{p\in\N_0}$ satisfying (lc), the condition  (dc) amounts to the condition $\log(M_p)=O(p^2)$, $p\to\infty$. Then, it suffices to apply \eqref{convolvedbelow}.
\item[(iii)] As seen in item (i), for every $t\ge 0$ we have
\begin{equation*}
2\omega_{\mathbb{M}}(t)=\omega_{\mathbb{M}\star\mathbb{M}}(t).
\end{equation*}
Since $\mathbb{M}$ satisfies (mg) if and only if there exists $H\ge 1$ such that
\begin{equation*}
2\omega_{\mathbb{M}}(t)\le\omega_{\mathbb{M}}(Ht)+H,\quad t\ge 0
\end{equation*}
(see \cite[Prop.~3.6]{komatsu}), it turns out that (mg) amounts to the fact that
$$
\omega_{\mathbb{M}\star\mathbb{M}}(t)\le \omega_{\mathbb{M}}(Ht)+H,\quad t\ge 0,
$$
for some $H\ge 1$, or in other words,
\begin{equation*}
h_{\mathbb{M}}(t)\le e^Hh_{\mathbb{M}\star\mathbb{M}}(Ht),\quad t>0.
\end{equation*}
\end{itemize}
\end{rema}

\subsection{Optimal flat functions and extension procedure}\label{subsectOptFlatConvolved}

Let $\mathbb{M}^1$ and $\mathbb{M}^2$ be weight sequences such that optimal flat functions $G_{\mathbb{M}^1}$ and $G_{\mathbb{M}^2}$ exist in the corresponding classes with uniform asymptotic expansion in a given sector $S$. Then, we claim that $G_{\mathbb{M}^1\star\mathbb{M}^2}:=G_{\mathbb{M}^1}\cdot G_{\mathbb{M}^2}$ is an optimal flat function (on the same sector $S$) in the class associated with the sequence $\mathbb{M}^1\star\mathbb{M}^2$.
Suppose $K_m$ and $J_m$, $m=1,2,3,4$, are the constants appearing in~\eqref{optimalflatleft} and \eqref{optimalflatright}
for $G_{\mathbb{M}^1}$ and $G_{\mathbb{M}^2}$, respectively. By \eqref{functionhforconvolved} we get that, for all $z\in S$,
$$|G_{\mathbb{M}^1}(z)\cdot G_{\mathbb{M}^2}(z)|\le K_3 h_{\mathbb{M}^1}(K_4|z|) J_3h_{\mathbb{M}^2}(J_4|z|)\le K_3J_3 h_{\mathbb{M}^1}(D|z|)h_{\mathbb{M}^2}(D|z|)= Ch_{\mathbb{M}^1\star\mathbb{M}^2}(D|z|),$$
with $C:=K_3J_3$ and $D:=\max\{K_4,J_4\}$, since each function $h_{\mathbb{M}}$ is nondecreasing. Similarly, for $x>0$ we can estimate
\begin{align*}
G_{\mathbb{M}^1}(x)\cdot G_{\mathbb{M}^2}(x)&\ge K_1 h_{\mathbb{M}^1}(K_2x) J_1h_{\mathbb{M}^2}(J_2x)\\
&\ge K_1J_1 h_{\mathbb{M}^1}(D_1x)h_{\mathbb{M}^2}(D_1x)= C_1h_{\mathbb{M}^1\star\mathbb{M}^2}(D_1x),
\end{align*}
with $C_1:=K_1J_1$ and $D_1:=\min\{K_2,J_2\}$, and the conclusion follows.

In case at least one of the sequences $\M^1$ and $\M^2$ satisfies (dc), $\mathbb{M}^1\star\mathbb{M}^2$ does so, and the extension operators from Theorem~\ref{tpral} will be available for the convolved sequence.

\subsection{Some examples}\label{convolveexamplesect}

Fix $q>1$ and $\sigma\in(1,2]$. Let us put $\L_{q,\sigma}:=\Ms\star\Ms$, $\L_{q,\sigma}=(L_p)_{p\in\N_0}$. It is not difficult to check that
$$
L_{2p}=q^{2p^{\sigma}},\quad L_{2p+1}=q^{{p^\sigma}+(p+1)^\sigma},
\qquad p\in\N_0.
$$
Observe that $2p^\sigma=2^{1-\sigma}(2p)^\sigma$, so that $L_{2p}$ equals the $2p$-th term of the sequence $\M_{q^{2^{1-\sigma}},\sigma}$. Regarding the odd terms, it is a consequence of Taylor's formula at $x=0$ for the functions of the form $x\mapsto (1+x)^{\alpha}$, $\alpha>0$, that
$$
p^{\sigma}+(p+1)^\sigma-2^{1-\sigma}(2p+1)^\sigma=O(p^{\sigma-2}),\quad p\to\infty.
$$
Since $\sigma\in(1,2]$, we deduce that $\L_{q,\sigma}$ is equivalent to $\M_{q^{2^{1-\sigma}},\sigma}$.

According to Subsection~\ref{subsectOptFlatConvolved}, an optimal flat function in the class associated with $\L_{q,\sigma}$ in, say, the sector $S_2$ is the function
$$
G(z):=G_2^{q,s}(z)G_2^{q,s}(z)= \exp\left(-2b_{q,s}\log^s\left(1+\frac{1}{z}\right)\right),\quad z\in S_2.
$$
It is not a surprise that, from the definition~\eqref{eq.def.b_qs} of $b_{q,s}$ and the relation between $\sigma$ and $s$, one obtains $b_{q^{2^{1-\sigma}},s}=2b_{q,s}$, and so
$G$ is precisely $G_2^{q^{2^{1-\sigma}},s}$, what agrees with the aforementioned equivalence of sequences.

If we consider instead $1<\sigma<2$ and $\mathbb{J}:=\M_{q,\sigma}\star\M_{q,2}$, $\mathbb{J}=(J_p)_{p\in\N_0}$, the computation of the terms $J_p$ is no longer possible in closed form, since their values depend for general $p$ on the position of $\sigma$ within the interval $(1,2)$. However, the previous subsection shows that, for $s$ associated with $\sigma$ as usual, the function
$$
G(z):=G_2^{q,s}(z)G_2^{q,2}(z)= \exp\left(-b_{q,s}\log^s\left(1+\frac{1}{z}\right)- b_{q,2}\log^2\left(1+\frac{1}{z}\right)\right),\quad z\in S_2,
$$
is an optimal flat function in the class associated with $\mathbb{J}$ in $S_2$. Note that $s$ is not equal to 2, hence the very aspect of the exponent in this function, and the fact that the restriction $G|_{(0,\infty)}$ is closely related to the function $h_{\mathbb{J}}$ (see Definition~\ref{optimalflatdef}), shows that $\mathbb{J}$ is not equivalent to any of the sequences $\Ms$. Since the sequence $\mathbb{J}$ does satisfy (dc), the extension procedure described in this paper is available for the classes associated with $\mathbb{J}$.

Observe that these examples of optimal flat functions can also be provided in general sectors $S_\ga$, $\ga>2$, by using the functions $G_{\ga}^{q,s}$ introduced in~\eqref{eq.defGgamma}.

\vskip.2cm
\noindent\textbf{Acknowledgements}: The first three authors are partially supported by the Spanish Ministry of Science and Innovation under the project PID2019-105621GB-I00. The fourth author is supported by FWF-Project P33417-N.

%
%
%

\vskip.5cm
\noindent\textbf{Affiliations}:\\
\noindent Javier~Jim\'{e}nez-Garrido:\\
Departamento de Matem\'aticas, Estad{\'\i}stica y Computaci\'on\\
Universidad de Cantabria\\
Avda. de los Castros, s/n, 39005 Santander, Spain\\
Instituto de Investigaci\'on en Matem\'aticas IMUVA, Universidad de Va\-lla\-do\-lid\\
ORCID: 0000-0003-3579-486X\\
E-mail: jesusjavier.jimenez@unican.es\\

\vskip.1cm
\noindent
Ignacio Miguel-Cantero:\\
Departamento de \'Algebra, An\'alisis Matem\'atico, Geometr{\'\i}a y Topolog{\'\i}a\\
Universidad de Va\-lla\-do\-lid\\
Facultad de Ciencias, Paseo de Bel\'en 7, 47011 Valladolid, Spain.\\
Instituto de Investigaci\'on en Matem\'aticas IMUVA\\
ORCID: 0000-0001-5270-0971\\
E-mail: ignacio.miguel@uva.es\\

\vskip.1cm
\noindent Javier~Sanz:\\
Departamento de \'Algebra, An\'alisis Matem\'atico, Geometr{\'\i}a y Topolog{\'\i}a\\
Universidad de Va\-lla\-do\-lid\\
Facultad de Ciencias, Paseo de Bel\'en 7, 47011 Valladolid, Spain.\\
Instituto de Investigaci\'on en Matem\'aticas IMUVA\\
ORCID: 0000-0001-7338-4971\\
E-mail: javier.sanz.gil@uva.es\\

\vskip.1cm
\noindent Gerhard~Schindl:\\
Fakult\"at f\"ur Mathematik, Universit\"at Wien,
Oskar-Morgenstern-Platz~1, A-1090 Wien, Austria.\\
ORCID: 0000-0003-2192-9110\\
E-mail: gerhard.schindl@univie.ac.at
\end{document}